\documentclass[12pt]{article}
\usepackage{amssymb}
\usepackage{amsmath}
\usepackage{array}
\usepackage[russian]{babel}
\usepackage[pdftex]{graphics}
\usepackage[X2,T2A]{fontenc}
\usepackage[cp1251]{inputenc}
\usepackage{mathrsfs}
\usepackage{longtable}
\usepackage{mathrsfs}
\usepackage{textcase}
\usepackage{setspace}
\usepackage{dsfont}
\usepackage{bbm}
\usepackage{pgf,pgfarrows,pgfnodes,pgfautomata,pgfheaps,pgfshade}
\usepackage{pdflscape}
\usepackage{longtable}

\DeclareMathOperator{\mmod}{mod}

\multlinegap=0em \DeclareFontFamily{T1}{msb}{}
\DeclareFontShape{T1}{msb}{m}{ol}{<5> <6> <7> <8> <9> gen * msbm
<10> <10.95> <12> <14.4> <17.28> <20.74> <24.88> msbm10}{}
\DeclareSymbolFont{AMSb}{T1}{msb}{m}{ol} \multlinegap=0em

\renewcommand{\S}{\mathhexbox278}

\begin{document}

\begin{center}
{\rmfamily\bfseries\Large A distribution related to Farey sequences - II}
\end{center}

\begin{center}
{\normalsize M.A.~Korolev\footnote{Steklov Mathematical Institute of Russian Academy of Sciences. 119991, Russia, Moscow, Gubkina str., 8. E-mail: \texttt{korolevma@mi-ras.ru}, \texttt{hardy\_ramanujan@mail.ru}}}
\end{center}

\fontsize{11}{13pt}\selectfont

We continue to study some arithmetical properties of Farey sequences by the method introduced by F.~Boca, C.~Cobeli and A.~Zaharescu (2001).
Let $\Phi_{Q}$ be the classical Farey sequence of order $Q$. Having the fixed integers $D\geqslant 2$ and $0\leqslant c_{0}\leqslant D-1$, we colour to the red the fractions in $\Phi_{Q}$ with denominators $\equiv c_{0} \; (\mmod D)$. Consider the gaps in $\Phi_{Q}$ with coloured endpoints, that do not contain the fractions $a/q$ with $q\equiv c_{0}\;(\mmod D)$ inside. The question is to find the limit proportions $\nu(r;D,c_{0})$ (as $Q\to +\infty$) of such gaps with precisely $r$ fractions inside in the whole set of the gaps under considering ($r = 0,1,2,3,\ldots$).

In fact, the expression for this proportion can be derived from the general result obtained by C.~Cobeli, M.~V\^{a}j\^{a}itu and A.~Zaharescu (2012). However, such formula expresses $\nu(r;D,c_{0})$ in the terms of areas of some polygons related to a special geometrical transform. In the present paper, we obtain explicit formulas for $\nu(r;D,c_{0})$ for the cases $3$ and $c_{0}=1,2$. Thus this and previous author's papers cover the case $D=3$.

\vspace{0.3cm}

\begin{center}
\textcolor{blue}{\textsc{Сontents}}
\end{center}

\textcolor{blue}{
\contentsline{subsection} {\numberline {\S\,1} Introduction}{\textcolor{blue}{1}}{}
\contentsline{subsection} {\numberline {\S\,2} Farey series, continuants and the sets $\mathcal{A}_{r}^{\circ}(D,c_{0},c_{1})$}{\textcolor{blue}{3}}{}
\contentsline{subsection} {\numberline {\S\,3} Basic properties of the polygons $\mathcal{T}(\mathbf{k})$}{\textcolor{blue}{7}}{}
\contentsline{subsection} {\numberline {\S\,4} Iterative construction of the sets $\mathcal{A}_{r}^{\circ}(3,1,2)$}{\textcolor{blue}{18}}{}
\contentsline{subsection} {\numberline {\S\,5} Iterative construction of the sets $\mathcal{A}_{r}^{\circ}(3,1,0)$}{\textcolor{blue}{22}}{}
\contentsline{subsection} {\numberline {\S\,6} Proof of the main assertion}{\textcolor{blue}{31}}{}
\contentsline{subsection} {\numberline {} Appendix I. Precise description of the polygons $\mathcal{T}(3,2,1,k)$, \\ $7\leqslant k\leqslant 12$}{\textcolor{blue}{35}}{}
\contentsline{subsection} {\numberline {} Appendix II. Precise description of the polygons $\mathcal{T}(3,2^{r-3},1,k)$, \\ $4r-10\leqslant k\leqslant 4r-4$}{\textcolor{blue}{36}}{}
\contentsline{subsection} {\numberline {} References}{\textcolor{blue}{39}}{}
}
\fontsize{12}{15pt}\selectfont
\vspace{1cm}

\section{Introduction}
\vspace{0.5cm}

In the present paper, we continue to study the arithmetic properties of Farey series started in \cite{Korolev_2023}, \cite{Korolev_2024}.\footnote{Preprint \cite{Korolev_2024} is a wide and modified version of the paper \cite{Korolev_2023}. That's why we use below some auxilliary assertions in the form stated in this preprint.} Let $\Phi_{Q}$ be Farey series of order $Q$, that is,
the set of all non-negative irreducible subunitary fractions whose denominators do not exceed $Q$, arranged in the ascending order. Further, let $r\geqslant 0$, $D\geqslant 2$, $0\leqslant c_{0}\leqslant D-1$ be integers. Let $N(Q;r,D,c_{0})$ denote the set of all tuples of $(r+2)$ consecutive fractions of series $\Phi_{Q}$ of the type
\begin{equation}\label{lab_01}
\frac{a_{0}}{q_{0}} < \frac{a_{1}}{q_{1}} < \ldots < \frac{a_{r}}{q_{r}} < \frac{a_{r+1}}{q_{r+1}},
\end{equation}
that satisfy to additional conditions
\begin{equation}\label{lab_02}
q_{0}, q_{r+1}\equiv c_{0} \pmod{D}, \quad q_{i}\not\equiv c_{0}\pmod{D}, \quad i=1,2,\ldots,r.
\end{equation}

Next, let $N(Q;D,c_{0})$ be the number of fractions of the above series whose denominators are congruent to $c_{0}$ modulo $D$. Consider the proportion
\[
\nu(Q;r,D,c_{0}) = \frac{N(Q;r,D,c_{0})}{N(Q;D,c_{0})}.
\]
F.~Boca, C.~Cobeli and A.~Zaharescu proposed in \cite{Boca_Cobeli_Zaharescu_2001} a new method for studying Farey series, which is based on the properties of so-called BCZ-transform $T$, which maps the <<Farey triangle>> $\mathcal{T} = \{(x,y)\,:\,0<x,y\leqslant 1, x+y>1$ into itself.

From the results of \cite{Cobeli_Zaharescu_2005}, \cite{Cobeli_Vajaitu_Zaharescu_2012} (see also \cite{Badziahin_Haynes_2011}), it follows that the limit
\[
\nu(r,D,c_{0}) = \lim_{Q\to +\infty}\nu(Q;r,D,c_{0})
\]
exists for any fixed $r$, $D$ and $c_{0}$.

Method of Boca, Cobeli and Zaharescu allows to express the limit proportions as the sums (infinite, generally speaking) of areas of some convex polygons defined in terms of BCZ-transform.
In general case, to write such proportion as the explicit function of $r, D, c_{0}$ seems to be quite difficult.

In the case $D=2$, $c_{0}=1$ and $1\leqslant r \leqslant 4$, such explicit formulas are given in \cite{Cobeli_Zaharescu_2005}, \cite{Cobeli_Vajaitu_Zaharescu_2012}.
These calculations were extended in \cite{Korolev_2024}. In the last paper, the case $D=3$, $c_{0} = 0$ was also considered. More precisely, the paper \cite{Korolev_2024} contains the asymptotic formulas for the quantities $\nu(Q;r,2,1)$, $r\geqslant 0$, and $\nu(Q;r,3,0)$, $r\geqslant 1$, valid in a quite wide range for $r$, namely, for all $r = o\bigl((Q/\ln{Q})^{1/3}\bigr)$.

Since the consecutive denominators of Farey fractions are coprime, they cannot be divisible by $3$ simultaneously. Hence, for any $Q$, the proportion $\nu(Q;r,3,0)$ vanishes identically for $r=0$
Thus, the case $D=3$, $c_{0}=0$ is studied completely.

In the present paper, we study the proportions $\nu(Q;r,3,c_{0})$ for $c_{0} = 1,2$. Our main result is following.
\vspace{0.5cm}

\textsc{Theorem.} \emph{Let $c_{0}$ be any of the numbers $1,2$. Then, for $Q\to +\infty$, $r\geqslant 0$, $r = o\bigl((Q/\ln{Q})^{1/3}\bigr)$ the following asymptotic formula holds:}
\[
\nu(Q;r,3,c_{0}) = \nu(r,3,c_{0}) + O\biggl(\frac{\ln{Q}}{Q}\biggr).
\]
\emph{Here the implied constant is absolute, and the quantities $\nu(r,3,c_{0})$ are defined as follows:}
\begin{align*}
& \nu(0;3,c_{0}) = \frac{1}{3} = 0.33333\,33333\ldots, \\[6pt]
& \nu(1;3,c_{0}) = \frac{2}{3}\biggl(\frac{\pi}{\sqrt{3}}-\ln{3}-\frac{1}{2}\biggr) = 0.14345\,80504\ldots,\\[6pt]
& \nu(2;3,c_{0}) = \frac{8}{3}\bigl(\ln{3}-1\bigr) = 0.26296\,61031\ldots, \\[6pt]
& \nu(3;3,c_{0}) = \frac{626}{105}-\frac{2\pi}{\sqrt{3}}-2\ln{3} = 0.13708\,14561\ldots, \\[6pt]
& \nu(4;3,c_{0}) = \frac{4}{3}\biggl(\frac{\pi}{\sqrt{3}}-\frac{23}{35}-\ln{3}\biggr) = 0.07739\,22912\ldots, \\[6pt]
& \nu(5;3,c_{0}) = \frac{4}{3}\ln{3} - \frac{193}{135} = 0.03518\,67552\ldots ,\\[6pt]
& \nu(r;3,c_{0}) = \frac{8}{3(2r-5)(2r-3)(2r-1)},\quad r\geqslant 6.
\end{align*}
Thus, the initial problem is completely solved for $D=2,3$. It is necessary to note that the case of the next prime modulus $D$, that is, $D=5$, seems much more complicated.
\vspace{0.5cm}

\section{Farey series, continuants and the sets $\boldsymbol{\mathcal{A}_{r}^{\circ}(D,c_{0},c_{1})}$}
\vspace{0.5cm}

As mentioned earlier, the general formula for the proportions $\nu(Q;r,D,c_{0})$ is given in \cite{Cobeli_Zaharescu_2005} (theorems 1, 2; see also \cite{Cobeli_Vajaitu_Zaharescu_2012}).
However, we need a modified version of such formula given in \cite{Korolev_2023}, which seems more appropriate for our purposes. Note that, for the case of fixed $r$, $D$ and $c_{0}$, the expressions for the limiting proportions $\nu(r,D,c_{0})$ that follows from lemma 23 of \cite{Korolev_2023} differs seemingly from the expressions that follows from theorems 1,2 of
\cite{Cobeli_Zaharescu_2005}. That's why we prove their equivalence in \cite[\S 9]{Korolev_2023}.

To formulate the corresponding result, we recall the basic properties of the BCZ-transformation $T: \mathcal{T}\to \mathcal{T}$, which is given by the formula
\begin{equation}\label{lab_03}
T(x,y) = \biggl(y,\biggl[\frac{1+x}{y}\biggr]y-x\biggr).
\end{equation}
The map $T$ is discontinuous, but it is bijective and area-preserving (see \cite[lemma 3]{Boca_Cobeli_Zaharescu_2001}). Moreover, if $a/q<a'/q'<a''/q''$ is the triple of consecutive Farey fractions of the series $\Phi_{Q}$, then (\ref{lab_03}) implies that
\[
T\biggl(\frac{q}{Q},\frac{q'}{Q}\biggr) = \biggl(\frac{q'}{Q},\frac{q''}{Q}\biggr).
\]
Further, given $k\geqslant 1$, denote by $\mathcal{T}(k)$ the set of all points of Farey triangle whose coordinates $(x,y)$ obey the condition $\bigl[(1+x)/y\bigr] = k$.
It is known that the closure $\mathcal{T}(k)$ is the triangle with vertices $(0,1)$, $(1,1)$ and $(\tfrac{1}{3},\tfrac{2}{3})$ for $k=1$, and the quadrangle with vertices
\[
A = \biggl(\frac{k}{k+2},\frac{2}{k+2}\biggr),\quad B = \biggl(1,\frac{2}{k}\biggr),\quad C =\biggl(1,\frac{2}{k+1}\biggr),\quad D = \biggl(\frac{k-1}{k+1},\frac{2}{k+1}\biggr),
\]
for $k\geqslant 2$ (see fig. 1 a,b). Further, given the tuple $\mathbf{k}_{r} = (k_{1},\ldots,k_{r})$ of integers $k_{j}\geqslant 1$, $j = 1,\ldots, r$, denote by
$\mathcal{T}(\mathbf{k}_{r})$ the set
\[
\mathcal{T}(k_{1})\bigcap T^{-1}\mathcal{T}(k_{2})\bigcap T^{-2}\mathcal{T}(k_{3})\bigcap\cdots \bigcap T^{-(r-1)}\mathcal{T}(k_{r}).
\]
If it is nonempty, then its closure is a convex polygon\footnote{In what follows, when mentioning the sets $\mathcal{T}(\mathbf{k}_{r})$, we omit the word <<closure>>.
Since the closure of any of such sets is obtained by the addition of some parts of its boundary, and the formulas for the proportions $\nu(Q;r,D,c_{0})$ involve the areas of such sets, this does not lead us to misunderstanding.}.

\begin{figure}[h]
\begin{minipage}[h]{1.0\linewidth}
\center{\includegraphics[width=0.5\linewidth]{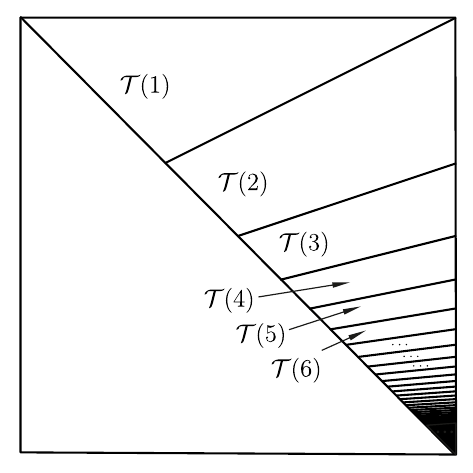}}\\
\end{minipage}
\begin{center}
Fig.~1а. \emph{The regions} $\mathcal{T}(k)$, $k\geqslant 1$.
\end{center}
\end{figure}

One can show that there exists a bijection between the pairs $q_{0}, q_{1}$ of consecutive denominators of Farey series $\Phi_{Q}$ and the primitive points of the dilated triangle $Q\!\cdot\!\mathcal{T}$ (see \cite[p.~562]{Boca_Cobeli_Zaharescu_2003}). Thus the computation of the number of the tuples satisfying to the conditions (\ref{lab_02}) reduces to the computation of the common number of primitive points lying in all nonempty regions $Q\!\cdot\!\mathcal{T}(\mathbf{k}_{r})$ corresponding to the tuples $\mathbf{k}_{r} = (k_{1},\ldots,k_{r})$ from some special set.

\begin{figure}[h]
\begin{minipage}[h]{1.0\linewidth}
\center{\includegraphics[width=0.9\linewidth]{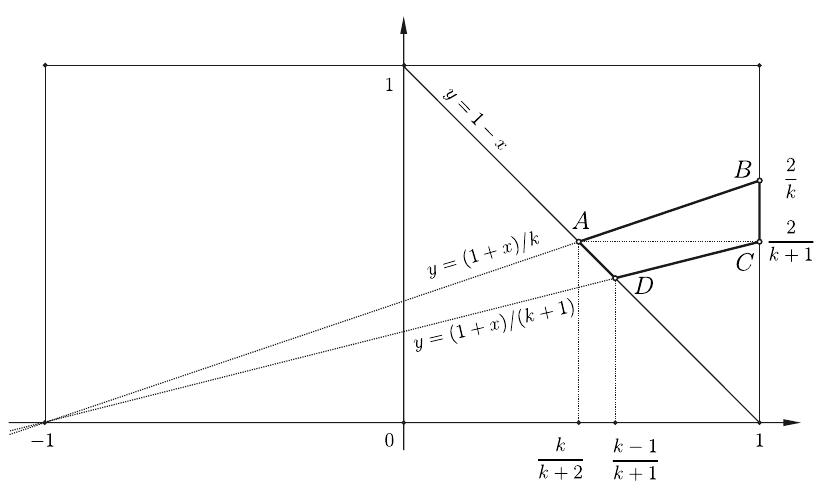}}\\
\end{minipage}
\begin{center}
Fig.~1b. \emph{Schematic drawing of the area $\mathcal{T}(k)$ for $k\geqslant 2$}.
\end{center}
\end{figure}

In turn, counting the number of primitive points in some <<good>> region on the coordinate plane ultimately comes down to calculation the area of such a region.
To calculate the area of the polygon $\mathcal{T}(k_{1},\ldots,k_{r})$, we use the fact that it is defined by a system of the inequalities of the type
\begin{equation}\label{lab_04}
0<x,y\leqslant 1,\quad x+y>1,\quad f_{i}(x;k_{1},\ldots,k_{i})<y\leqslant g_{i}(x;k_{1},\ldots,k_{i}),\quad i = 1,2,\ldots,r,
\end{equation}
where
\begin{equation}\label{lab_05}
f_{i}(x;k_{1},\ldots,k_{i}) = \frac{1+x\,\mathbb{K}_{i-1}(k_{2},\ldots,k_{i}+1)}{\mathbb{K}_{i}(k_{1},\ldots,k_{i}+1)},\quad
g_{i}(x;k_{1},\ldots,k_{i}) = \frac{1+x\,\mathbb{K}_{i-1}(k_{2},\ldots,k_{i})}{\mathbb{K}_{i}(k_{1},\ldots,k_{i})}.
\end{equation}
The arguments given in \cite{Korolev_2024} and based on the method of \cite{Boca_Cobeli_Zaharescu_2001}-\cite{Cobeli_Vajaitu_Zaharescu_2012} lead to the following expression for the
proportion $\nu(Q;r,D,c_{0})$:
\begin{equation}\label{lab_06}
\nu(Q;r,D,c_{0}) = \frac{2}{D}\frac{\Delta}{\varphi(\Delta)}\bigl(\mathfrak{c}_{r} + O(R)\bigr),
\end{equation}
where
\begin{align*}
& \Delta = \text{GCD}(c_{0},D),\quad \mathfrak{c}_{r} = \sum\limits_{\substack{0\leqslant c_{1}\leqslant D-1 \\ c_{1}\ne c_{0} \\ \text{GCD}(c_{1},\Delta)=1}}\sum\limits_{\mathbf{k}\in\,\mathcal{A}_{r}^{\circ}(D,c_{0},c_{1})}|\mathcal{T}(\mathbf{k})|, \\
& R = \frac{D^{2}}{Q}\biggl(1+\frac{\ln{Q}}{D}\biggr)\Sigma_{0} + \bigl(D^{3}+Dr\ln{Q}+r(\ln{Q})^{2}\bigr)\frac{\mathbbm{1}_{\{r\leqslant 2D\}}}{Q^{2}}+\mathfrak{c}_{r}\,\frac{D\Delta}{\varphi(\Delta)}\,\frac{\ln{Q}}{Q}\\[6pt]
& \Sigma_{0} = \sum\limits_{\substack{0\leqslant c_{1}\leqslant D-1 \\ c_{1}\ne c_{0} \\ \text{GCD}(c_{1},\Delta)=1}}\sum\limits_{\substack{\mathbf{k}\in\, \mathcal{A}_{r}^{\circ}(D,c_{0},c_{1}) \\ \|\mathbf{k}\|\leqslant 4r+1}}1.
\end{align*}
Here $\mathcal{A}_{r}^{\circ}(D,c_{0},c_{1})$ denotes the set of all tuples $\mathbf{k} = (k_{1},\ldots,k_{r})$ satisfying the conditions
\begin{equation}\label{lab_07}
\begin{cases}
c_{1}\,\mathbb{K}_{i}(k_{1},\ldots,k_{i})-c_{0}\,\mathbb{K}_{i-1}(k_{2},\ldots,k_{i})\not\equiv c_{0}\pmod{D},\quad i = 1,2,\ldots, r-1,\\
c_{1}\,\mathbb{K}_{r}(k_{1},\ldots,k_{r})-c_{0}\,\mathbb{K}_{r-1}(k_{2},\ldots,k_{r})\equiv c_{0}\pmod{D}
\end{cases}
\end{equation}
and such that the region $\mathcal{T}(\mathbf{k})$ is nonempty. Here and below, the symbol $\|\mathbf{k}\|$ means the maximal component of the tuple $\mathbf{k}$:
$\|\mathbf{k}\| = \max{\{k_{1},\ldots,k_{r}\}}$. The quantity $\mathbbm{1}_{A}$ equals to one if the condition $A$ holds, and equals to zero otherwise.
 The constant in $O$-symbol is absolute, and the formula (\ref{lab_06}) is valid for any $r$ and $D$ such that $1\leqslant r\leqslant Q/3$, $D\geqslant 2$, $D = o\bigl(\sqrt{Q\mathstrut}\bigr)$, and for any $c_{0}$, $0\leqslant c_{0}\leqslant D-1$.

In particular, if the difference $D$ of the progression is fixed (do not depend on $Q$), the bound for $R$ from (\ref{lab_06}) reduces to
\begin{equation}\label{lab_08}
R \ll_{D} \bigl(\mathfrak{c}_{r}+\Sigma_{0}\bigr)\frac{\ln{Q}}{Q}.
\end{equation}

In \cite{Korolev_2024}, we give an explicit description of the sets $\mathcal{A}_{r}^{\circ}(3;0,c_{1})$, $r \geqslant 1$, $c_{1} = 1,2$.
It is the main ingredient that allows to obtain the asymptotic formulas for the proportions $\nu(Q;r,3,0)$. In the present paper, we give an explicit description of each set
$\mathcal{A}_{r}^{\circ}(3;c_{0},c_{1})$, $c_{0} = 1,2$, $0\leqslant c_{1}\leqslant 2$, $c_{1}\ne c_{0}$, and the calculation of the areas of the polygons $\mathcal{T}(\mathbf{k})$
that correspond to all tuples $\mathbf{k} \in \mathcal{A}_{r}^{\circ}(3;c_{0},c_{1})$.

First we have to transform the formula (\ref{lab_06}) to more appropriate form.
\vspace{0.5cm}

\textsc{Lemma 1.} \emph{In each case $c_{0}=1$, $c_{0} = 2$, the relation}
\[
\nu(Q;r,3,c_{0}) = \frac{2}{3}\,\mathfrak{c}_{r} + O(R),\quad\textit{where}\quad
\mathfrak{c}_{r} = \biggl(\,\sum\limits_{\mathbf{k}\in\,\mathcal{A}_{r}^{\circ}(3,1,0)} + \sum\limits_{\mathbf{k}\in\,\mathcal{A}_{r}^{\circ}(3,1,2)}\;\biggr)|\mathcal{T}(\mathbf{k})|,
\]
\emph{holds for any $r$, $1\leqslant r\leqslant Q/3$. Here the term $R$ obeys the estimate} (\ref{lab_08}), \emph{and the implied constant is absolute}.
\vspace{0.5cm}

\textsc{Proof.} In each case $c_{0}=1$, $c_{0} = 2$ we have $\Delta =$ $\text{GCD}(c_{0},3)$ $=1$, so (\ref{lab_06}) takes the form
\[
\nu(Q;r,3,c_{0}) = \frac{2}{3}\,\mathfrak{c}_{r} + O(R).
\]
Further,
\[
\mathfrak{c}_{r} = \biggl(\,\sum\limits_{\mathbf{k}\in\,\mathcal{A}_{r}^{\circ}(3,1,0)} + \sum\limits_{\mathbf{k}\in\,\mathcal{A}_{r}^{\circ}(3,1,2)}\;\biggr)|\mathcal{T}(\mathbf{k})|
\]
in the case $c_{0}=1$, and
\[
\mathfrak{c}_{r} = \biggl(\,\sum\limits_{\mathbf{k}\in\,\mathcal{A}_{r}^{\circ}(3,2,0)} + \sum\limits_{\mathbf{k}\in\,\mathcal{A}_{r}^{\circ}(3,2,1)}\;\biggr)|\mathcal{T}(\mathbf{k})|
\]
in the case $c_{0}=2$.

Now we show that, for any $r\geqslant 1$, both the sets $\mathcal{A}_{r}^{\circ}(3,1,2)$ and $\mathcal{A}_{r}^{\circ}(3,2,1)$ consist of the same tuples, and the same is true for the sets
$\mathcal{A}_{r}^{\circ}(3,1,0)$ and $\mathcal{A}_{r}^{\circ}(3,2,0)$.

In the case $D = 3, c_{0}=1, c_{1} = 2$, the conditions (\ref{lab_07}) for the tuple $\mathbf{k} = (k_{1},\ldots, k_{r})$ to belong to $\mathcal{A}_{r}^{\circ}(3,1,2)$ have the form
\begin{equation*}
\begin{cases}
2\mathbb{K}_{i}(k_{1},\ldots,k_{i})-\mathbb{K}_{i-1}(k_{2},\ldots,k_{i})\not\equiv 1\pmod{3},\quad i = 1,2,\ldots, r-1,\\
2\mathbb{K}_{r}(k_{1},\ldots,k_{r})-\mathbb{K}_{r-1}(k_{2},\ldots,k_{r})\equiv 1\pmod{3},
\end{cases}
\end{equation*}
or, that is the same,
\begin{equation}\label{lab_09}
\begin{cases}
\mathbb{K}_{i}(k_{1},\ldots,k_{i})+\mathbb{K}_{i-1}(k_{2},\ldots,k_{i})\not\equiv 2\pmod{3},\quad i = 1,2,\ldots, r-1,\\
\mathbb{K}_{r}(k_{1},\ldots,k_{r})+\mathbb{K}_{r-1}(k_{2},\ldots,k_{r})\equiv 2\pmod{3}.
\end{cases}
\end{equation}
Further, the conditions for a tuple to belong to the set $\mathcal{A}_{r}^{\circ}(3,2,1)$ have the form
\begin{equation*}
\begin{cases}
\mathbb{K}_{i}(k_{1},\ldots,k_{i})-2\mathbb{K}_{i-1}(k_{2},\ldots,k_{i})\not\equiv 2\pmod{3},\quad i = 1,2,\ldots, r-1,\\
\mathbb{K}_{r}(k_{1},\ldots,k_{r})-2\mathbb{K}_{r-1}(k_{2},\ldots,k_{r})\equiv 2\pmod{3}
\end{cases}
\end{equation*}
and are obviously equivalent to (\ref{lab_09}).

Also, it is easy to check that the conditions for the tuple $\mathbf{k} = (k_{1},\ldots, k_{r})$ to belong to the sets
$\mathcal{A}_{r}^{\circ}(3,1,0)$ и $\mathcal{A}_{r}^{\circ}(3,2,0)$ reduce to the same form, namely
\begin{equation}\label{lab_10}
\begin{cases}
\mathbb{K}_{i-1}(k_{2},\ldots,k_{i})\not\equiv 2\pmod{3},\quad i = 1,2,\ldots, r-1,\\
\mathbb{K}_{r-1}(k_{2},\ldots,k_{r})\equiv 2\pmod{3}.
\end{cases}
\end{equation}
Lemma is proved.
\vspace{0.5cm}

\textsc{Remark.} In the case $r=1$, the condition (\ref{lab_10}) has the form $\mathbb{K}_{0}\equiv 2\pmod{3}$ and fails in view of the relation $\mathbb{K}_{0}=1$.
This means that the sets $\mathcal{A}_{1}^{\circ}(3,1,0)$ and $\mathcal{A}_{1}^{\circ}(3,2,0)$ are empty. In the case $r\geqslant 2$, it is sufficient to check the first of the conditions
(\ref{lab_10}) for $i\geqslant 2$ only, because for $i=1$ it is fulfilled automatically.
\vspace{0.5cm}

\section{Basic properties of the polygons $\boldsymbol{\mathcal{T}(\mathbf{k})}$}
\vspace{0.5cm}

This section contains the auxilliary assertions concerning the basic properties of the regions $\mathcal{T}(k_{1},\ldots,k_{r})$.
Below, checking the non-emptiness of such regions (see basic lemmas 13, 14) we will use repeatedly the following simple argument.
Suppose that $r\geqslant 3$ and let $2\leqslant j\leqslant r-1$. Since
\begin{multline*}
\mathcal{T}(k_{1},\ldots,k_{r}) = \mathcal{T}(k_{1})\bigcap T^{-1}\mathcal{T}(k_{2})\bigcap \cdots \bigcap T^{-(j-2)}\mathcal{T}(k_{j-1})\bigcap \\
\bigcap T^{-(j-1)}\bigl\{\mathcal{T}(k_{j})\bigcap  T^{-1}\mathcal{T}(k_{j+1})\ldots \bigcap T^{-(r-j-1)}\mathcal{T}(k_{r-j})\bigr\} = \\
= \mathcal{T}(k_{1},\ldots,k_{j-1})\bigcap T^{-(j-1)}\mathcal{T}(k_{j},\ldots,k_{r}),
\end{multline*}
then the non-emptiness of both regions
\[
\mathcal{T}(k_{1},\ldots,k_{j-1}),\quad \mathcal{T}(k_{j},\ldots,k_{r})
\]
is the necessary condition for non-emptiness of the initial region $\mathcal{T}(k_{1},\ldots,k_{r})$.
The repeated application of this argument leads to a test of non-emptiness of sets that correspond to tuples of small length: $r = 2,3$.
The complete classification of non-empty regions for such $r$ is given by the next two assertions from \cite{Korolev_2024} (lemmas 10, 12).
\vspace{0.5cm}

\textsc{Lemma 2.} \emph{The region $\mathcal{T}(k,m)$ is empty iff the pair $(k,m)$ satisfies to one of the following conditions:}\\ [6pt]
\begin{tabular}{l p{0.5\linewidth}}
(a) $m=1$, $k=1$;  & (b) $m=2$, $k\geqslant 5$;  \\[6pt]
(c) $m=3,4$, $k\geqslant 3$; & (d) $m\geqslant 5$, $k\geqslant 2$.
\end{tabular}
\vspace{0.5cm}

\textsc{Lemma 3.} \emph{The region $\mathcal{T}(k,m,n)$ is empty iff the triple $(k,m,n)$ satisfies to at least one of the following conditions:} \\ [6pt]
\begin{tabular}{l p{0.5\linewidth}}
(1) \emph{the pair $(k,m)$ or $(m,n)$ satisfies to the conditions of lemma 2}; & \\ [6pt]
\end{tabular}
\begin{tabular}{l p{0.5\linewidth}}
(2a) $m=1$, $k=2$, $1\leqslant n\leqslant 5$;  & (2b) $m=1$, $k=3$, $n\ne 4,5,6,7,8$;\\[6pt]
(2c) $m=1$, $k=4$, $n\ne 3,4,5$;  & (2d) $m=1$, $k=5$, $n\ne 3,4$; \\[6pt]
(2e) $m=1$, $6\leqslant k\leqslant 8$, $n\geqslant 4$; & (2f) $m=1$, $k\geqslant 9$, $n\geqslant 3$; \\[6pt]
(3a) $m=2$, $k=1$, $n\ne 2,3,4$; & (3b) $m=2$, $k=2$, $n\geqslant 4$; \\[6pt]
(3c) $m=2$, $k=3$, $n\geqslant 3$; & (3d) $m=2$, $k=4$, $n\geqslant 2$; \\[6pt]
(4a) $m=3$, $k=1$, $n\geqslant 3$; & (4b) $m=3$, $k=2$, $n\geqslant 3$; \\[6pt]
(5a) $m=4$, $k=1$, $n\geqslant 3$; & (5b) $m=4$, $k=2$, $n\geqslant 2$; \\[6pt]
(6)  $m\geqslant 5$, $k=1$, $n\geqslant 2$.
\end{tabular}
\vspace{0.5cm}

In some cases, one can reduce the searching of non-empty regions $\mathcal{T}(k_{1},\ldots,k_{r})$ by using two following assertions which are, correspondingly, the corollary 1 of lemma 11 and p. a of lemma 2 from \cite{Korolev_2024}.
\vspace{0.5cm}

\textsc{Lemma 4.} \emph{For any positive integers $k_{1},\ldots,k_{r}$, one has}
\[
|\mathcal{T}(k_{r},\ldots,k_{1})| = |\mathcal{T}(k_{1},\ldots,k_{r})|.
\]

\textsc{Lemma 5.} \emph{For non-empty region $\mathcal{T}(k_{1},\ldots,k_{r})$, one has $\mathbb{K}_{r}(k_{1},\ldots,k_{r})\geqslant 1$.}
\vspace{0.3cm}

The following assertion shows that if the tuple $(k_{1},\ldots,k_{r})$ contains a quite large component $k_{j}$, then such tuple has a very special form.
\vspace{0.3cm}

\textsc{Lemma 6.} \emph{Let} $r\geqslant 1$, $\max{\{k_{1},\ldots,k_{r}\}} = k_{j}\geqslant 4r+2$, \emph{and suppose that} $\mathcal{T}(k_{1},\ldots,k_{r})$ \emph{is non-empty. Then the tuple} $(k_{1},\ldots,k_{r})$ \emph{has the form:}
\begin{equation}\label{lab_11}
k_{i} = 1\quad\textit{при} \quad |i-j|=1,\quad k_{i} = 2\quad\textit{при} \quad |i-j|\geqslant 2.
\end{equation}

For the proof, see: \cite[lemma 3.4 (\textit{i})]{Boca_Cobeli_Zaharescu_2003} (see also the 2nd passage on p.~570)\footnote{In \cite{Boca_Cobeli_Zaharescu_2003}, the notation $\mathcal{T}_k$ is used instead of $\mathcal{T}(k)$.}.

The following lemma is useful for the evaluation of some infinite series that arise in the process of calculation of limit proportions $\nu(r;3,c_{0})$, $c_{0} = 1,2$.
\vspace{0.3cm}

\textsc{Lemma 7.} \emph{The following relations hold true:} \\[6pt]
\begin{tabular}{l p{0.5\linewidth}}
(a) $|\mathcal{T}(k,1)| = |\mathcal{T}(k)|$ \emph{for} $k\geqslant 5$;  & (b) $|\mathcal{T}(k,1,2)| = |\mathcal{T}(k)|$ \emph{for} $k\geqslant 9$; \\ [6pt]
(c) $|\mathcal{T}(2,1,k)| = |\mathcal{T}(k)|$ \emph{for} $k\geqslant 9$; & (d) $|\mathcal{T}(2,1,k,1)| = |\mathcal{T}(k)|$ \emph{for} $k\geqslant 9$; \\[6pt]
(e) $|\mathcal{T}(2,1,k,1,2)| = |\mathcal{T}(k)|$ \emph{for} $k\geqslant 9$; & (f) $|\mathcal{T}(1,k,1)| = |\mathcal{T}(k)|$ \emph{for} $k\geqslant 6$.\\[6pt]
\end{tabular}

This is lemma 12 from \cite{Korolev_2024}.
\vspace{0.5cm}

Finally, to describe the sets $\mathcal{A}_{r}^{\circ}(3,c_{0},c_{1})$, we need several lemmas concerning the emptiness (non-emptiness) of some regions $\mathcal{T}(\mathbf{k})$ that correspond to the tuples $\mathbf{k}_{r}  = (k_{1},\ldots,k_{r})$ of special type. Its proofs are similar and follow the same lines. Having a system of the inequalities describing the polygon
$\mathcal{T}(k_{1},\ldots,k_{r-1})$, we use an explicit formulas for the linear functions $f(x) = f_{r}(x;\mathbf{k}_{r})$, $g(x) = g_{r}(x;\mathbf{k}_{r})$ (see (\ref{lab_05}))) to find the values of $k_{r}$ such that the intersection of $\mathcal{T}(k_{1},\ldots,k_{r-1})$ with the region defined by the inequalities $f(x)<y\leqslant g(x)$, in non-empty.
\vspace{0.3cm}

\textsc{Lemma 8.} \emph{The region $\mathcal{T}(k,1,2,4)$ is non-empty for $6\leqslant k\leqslant 8$ only.}
\vspace{0.3cm}

\textsc{Proof.} By lemma 4, both the regions $\mathcal{T}(k,1,2,4)$ and $\mathcal{T}(4,2,1,k)$ are simultaneously empty or non-empty. By lemma 3 (2a), the necessary condition of non-emptiness of the last region, that is, the non-emptiness of $\mathcal{T}(2,1,k)$, implies the inequality $k\geqslant 6$. Further, by direct calculation, one can check that the region $\mathcal{T}(4,2,1)$ is described by the system
\begin{equation}\label{lab_12}
\frac{2}{3}<x\leqslant 1,\quad \frac{1+x}{5}<y\leqslant \frac{1+2x}{7}
\end{equation}
and appears to be the triangle with vertices
\[
A = \biggl(\frac{2}{3},\frac{1}{3}\biggr),\quad B = \biggl(1,\frac{3}{7}\biggr),\quad C = \biggl(1,\frac{2}{5}\biggr).
\]
When passing to the region $\mathcal{T}(4,2,1,k)$, the inequalities $f(x)<y\leqslant g(x)$,
\begin{multline*}
f(x) = f_{4}(x;4,2,1,k) = \frac{1+x\,\mathbb{K}_{3}(2,1,k+1)}{\mathbb{K}_{4}(4,2,1,k+1)} = \frac{1+x(k-1)}{3k-4},\\
g(x) = g_{4}(x;4,2,1,k) = \frac{1+x\,\mathbb{K}_{3}(2,1,k)}{\mathbb{K}_{4}(4,2,1,k)} = \frac{1+x(k-2)}{3k-7},
\end{multline*}
are added to the system (\ref{lab_12}).

The solutions of the inequality $g(x)\leqslant (1+x)/5$ form a ray $x\leqslant x_{k}$, where
\[
x_{k} = \frac{3k-12}{2k-3}.
\]
For any $k\geqslant 9$, such ray contains the segment $2/3\leqslant x\leqslant 1$, which is the projection of $\triangle ABC$ onto horizontal axis. Therefore, if $k\geqslant 9$
then the triangle $ABC$ lies over the line $y = g(x)$, so the region $\mathcal{T}(4,2,1,k)$ is empty. The non-emptiness of $\mathcal{T}(4,2,1,k)$ for $6\leqslant k\leqslant 8$ is checked by direct calculation (see. Fig. 2).

\begin{figure}[h]
\begin{minipage}[h]{1.0\linewidth}
\center{\includegraphics[width=0.7\linewidth]{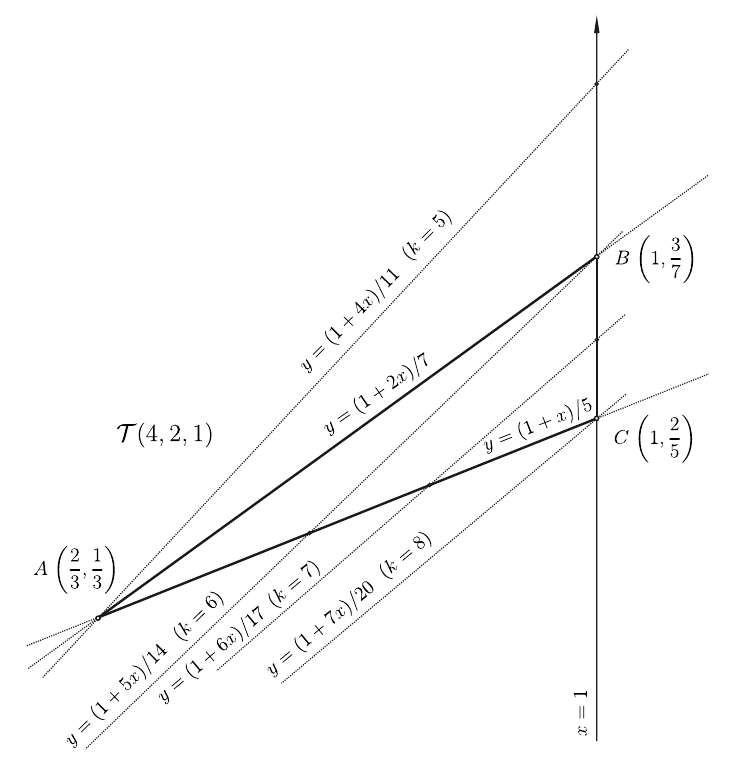}}\\
\end{minipage}
\begin{center}
Fig.~2. \emph{Schematic drawing of the polygon $\mathcal{T}(4,2,1)$ and the lines $y = f(x)$ corresponding to $5\leqslant k\leqslant 8$}.
\end{center}
\end{figure}

\textsc{Lemma 9.} \emph{Suppose that $n\geqslant 2$. Then}\\
(a) \emph{the region $\mathcal{T}(3,2^{n})$ is described by the system}
\begin{equation}\label{lab_13}
\frac{2n-1}{2n+1}<x\leqslant 1,\quad \frac{1+x}{4}<y\leqslant \frac{1+(n+1)x}{2n+3}
\end{equation}
\emph{and is the triangle with vertices}
\[
A = \biggl(\frac{2n-1}{2n+1},\frac{n}{2n+1}\biggr),\quad B = \biggl(1,\frac{n+2}{2n+3}\biggr),\quad C = \biggl(1,\frac{1}{2}\biggr);
\] \\
(b) \emph{the region $\mathcal{T}(3,2^{n},k)$ is empty for $k \geqslant 3$ and is non-empty for $1\leqslant k\leqslant 2$.}
\vspace{0.3cm}

\textsc{Proof.} We prove (a) by induction with respect to $n$ and, simultaneously, check the validity of (b). For $n=2$, the assertion (a) is checked by direct calculation. Suppose that it holds for all $n$, $2\leqslant n\leqslant m$, and prove that (b) is valid  for $n = m$, and (a) is valid for $n = m+1$.

Indeed, when passing from $\mathcal{T}(3,2^{m})$ to the region $\mathcal{T}(3,2^{m},k)$, we need to add the inequalities $f(x)<y\leqslant g(x)$ to the system (\ref{lab_13}); here
\begin{multline*}
f(x) = f_{m+2}(x;3,2^{m},k) = \frac{1+x\,\mathbb{K}_{m+1}(2^{m},1,k+1)}{\mathbb{K}_{m+2}(3,2^{m},k+1)} = \frac{1+x(km+k+1)}{2km+3k+2},\\
g(x) = g_{m+2}(x;3,2^{m},k) = \frac{1+x\,\mathbb{K}_{m+1}(2^{m},1,k)}{\mathbb{K}_{m+2}(3,2^{m},k)} = \frac{1+x(km+k-m)}{2km+3k-2m-1}
\end{multline*}
(here and below, we use an explicit formulas for the continuants given by lemmas 19, 20 of \cite{Korolev_2024}).
It is easy to check that the inequality $f(x)<g(x)$ holds for $x>-(2m+3)$ and, in particular, it is satisfied by any $x$, $0\leqslant x\leqslant 1$.
Further, the solutions of the inequality $g(x)<(1+x)/4$ form the ray $x\leqslant x_{k}$, where
\[
x_{k} = \frac{(2m+1)(k-1)+2+2(k-3)}{(2m+1)(k-1)+2}.
\]
In particular, if $k\geqslant 3$ then, for any $x$, $0\leqslant x\leqslant 1$, we have $f(x)<g(x)<(1+x)/4$. Hence, the region $\mathcal{T}(3,2^{m},k)$ is empty for such $k$.

The values $k = 1,2$ lead the pairs
\begin{multline*}
f(x) = \frac{1+(m+2)x}{2m+5},\quad g(x) = \frac{1+x}{2}\quad \text{and}\\
f(x) = \frac{1+(2m+3)x}{4m+8},\quad g(x) = \frac{1+(m+2)x}{2m+5},
\end{multline*}
correspondingly. It is easy to check that $\triangle ABC$ lies entirely under the line $y = (1+x)/2$ and lies entirely over the line $y = (1+(2m+3)x)/(4m+8)$.
At the same time, the line
\[
y = \frac{1+(m+2)x}{2m+5}
\]
meets with the segments $AB$ and $AC$ at the points
\[
A_{1} = \biggl(\frac{2m+1}{2m+3},\frac{m+1}{2m+3}\biggr),\quad B_{1} = \biggl(1,\frac{m+3}{2m+5}\biggr)
\]
and thus cuts $\triangle ABC$ into two parts. The upper one corresponds to the value $k=1$ while the second one satisfies to $k=2$.
Thus, both the regions $\mathcal{T}(3,2^{m},1)$, $\mathcal{T}(3,2^{m+1})$ are non-empty and described by the systems
 \begin{equation}\label{lab_14}
\begin{cases}
\displaystyle \frac{2m-1}{2m+1} < x\leqslant \frac{2m+1}{2m+3},\quad \frac{1+x}{4} < y\leqslant \frac{1+(m+1)x}{2m+3},\\[12pt]
\displaystyle \frac{2m+1}{2m+3} < x\leqslant 1,\quad \frac{1+(m+2)x}{2m+5} < y\leqslant \frac{1+(m+1)x}{2m+3}
\end{cases}
\end{equation}
and
\[
\frac{2m+1}{2m+3} < x\leqslant 1,\quad \frac{1+(m+2)x}{2m+5} < y\leqslant \frac{1+(m+2)x}{2m+5}
\]
(see. fig. 3). Lemma is proved.

\begin{figure}[h]
\begin{minipage}[h]{1.0\linewidth}
\center{\includegraphics[width=0.7\linewidth]{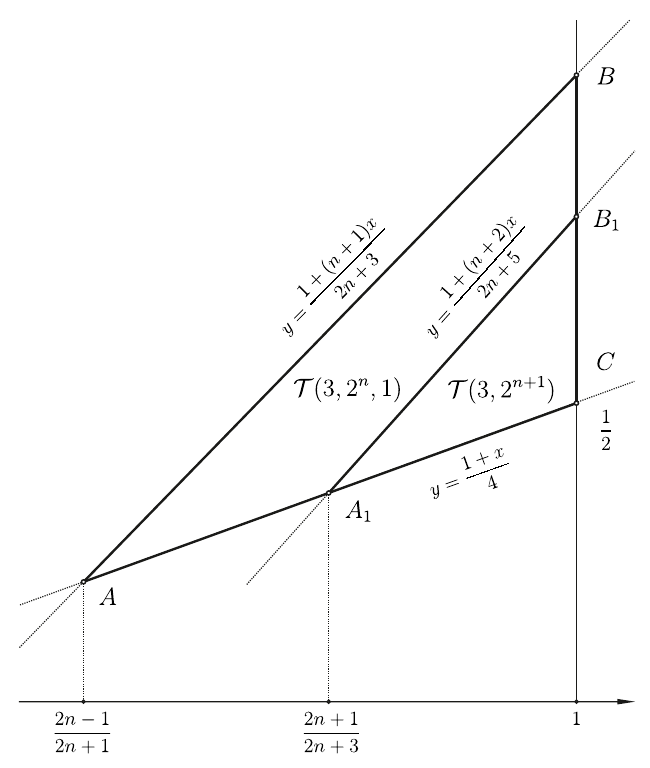}}\\
\end{minipage}
\begin{center}
Fig.~3. \emph{Schematic drawing of the polygons $\mathcal{T}(3,2^{n},1)$ and $\mathcal{T}(3,2^{n+1})$}.
\end{center}
\end{figure}

\textsc{Lemma 10.} \emph{Let $n\geqslant 1$. The region $\mathcal{T}(k,1,2^{n})$ is non-empty iff $k\geqslant 4n+2$.}
\vspace{0.3cm}

\textsc{Proof.} By lemma 4, both the regions $\mathcal{T}(k,1,2^{n})$ and $\mathcal{T}(2^{n},1,k)$ are simultaneously empty or non-empty.
In $n=1$ then the desired assertion follows from lemma 3 (2a). Hence, we assume further that $n\geqslant 2$.

Since $\mathbb{K}_{n+2}(2^{n},1,k) = k-n-1\leqslant 0$ for $1\leqslant k\leqslant n+1$, then lemma 5 implies that $\mathcal{T}(2^{n},1,k)$ is empty for such $k$.
So, we may assume that $k\geqslant n+2$.

Arguing as in the proof of lemma 9, one can show that the region $\mathcal{T}(2^{n},1)$ is described by the system
\begin{equation}\label{lab_15}
\begin{cases}
\displaystyle \frac{n}{2n+1} < x\leqslant \frac{n+1}{2n+3},\quad 1-x < y\leqslant \frac{1+nx}{n+1},\\[12pt]
\displaystyle \frac{n+1}{2n+3} < x\leqslant 1,\quad \frac{1+(n+1)x}{n+2} < y\leqslant \frac{1+nx}{n+1}
\end{cases}
\end{equation}
and is the triangle with vertices
\[
A = \biggl(\frac{n}{2n+1},\frac{n+1}{2n+1}\biggr),\quad B = \bigl(1,1\bigr),\quad C = \biggl(\frac{n+1}{2n+3},\frac{n+2}{2n+3}\biggr)
\]
(see. fig. 4а).

\begin{figure}[h]
\begin{minipage}[h]{1.0\linewidth}
\center{\includegraphics[width=0.7\linewidth]{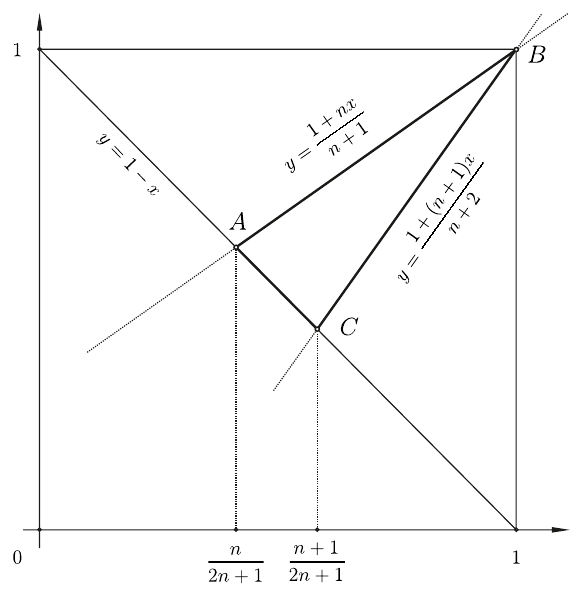}}\\
\end{minipage}
\begin{center}
Fig.~4а. \emph{Schematic drawing of the polygon $\mathcal{T}(2^{n},1)$}.
\end{center}
\end{figure}

When passing to the region $\mathcal{T}(2^{n},1,k)$, the inequalities $F_{k}(x)<y\leqslant G_{k}(x)$,
\begin{multline*}
F_{k}(x) = f_{n+2}(x;2^{n},1,k) = \frac{1+x\,\mathbb{K}_{n+1}(2^{n-1},1,k+1)}{\mathbb{K}_{n+2}(2^{n},1,k+1)} = \frac{1+x(k-n+1)}{k-n}, \\
G_{k}(x) = g_{n+2}(x;2^{n},1,k) = \frac{1+x\,\mathbb{K}_{n+1}(2^{n-1},1,k)}{\mathbb{K}_{n+2}(2^{n},1,k)} = \frac{1+x(k-n)}{k-n-1}
\end{multline*}
are added to the system (\ref{lab_15}). However, for any $k$, $n+2\leqslant k\leqslant 4n+1$, the line $y = F_{k}(x)$ lies over $\triangle ABC$. Indeed, the solutions of the inequality $(1+nx)/(n+1)\leqslant F_{k}(x)$ form the ray $x\geqslant x_{k}$, where
\[
x_{k} = \frac{k-2n-1}{k+1}.
\]
If the condition $x_{k}\leqslant n/(2n+1)$ holds, then it contains the projection of $\triangle ABC$ onto the horizontal axis, that is, the segment
\[
\frac{n}{2n+1}\leqslant x\leqslant 1,
\]
But such condition is valid for all $k\leqslant 4n+1$. Hence, the region $\mathcal{T}(2^{n},1,k)$ is empty for $n+2\leqslant k\leqslant 4n+1$.
It is easy to check that
\[
1< 1+\frac{2}{k-n} = F_{k}(1) < G_{k}(1) = 1 + \frac{2}{k-n-1}.
\]
Further, the $x$-coordinates of the points of intersection of the lines $y = F_{k}(x)$, $y = G_{k}(x)$ with the line $y = 1-x$ are equal to
\[
X_{k} = \frac{k-n-1}{2(k-n)+1},\quad X_{k-1} = \frac{k-n-2}{2(k-n)-1},
\]
correspondingly, and satisfy the inequalities
\[
\frac{n}{2n+1}<X_{k-1}<X_{k}<1,
\]
for any $k\geqslant 4n+2$. Hence, both the lines cut the $\triangle ABC$ (see fig. 4b). Finally, it is easy to check that the inequality
$F_{k}(x)\leqslant G_{k}(x)$ holds for any $x\geqslant -1$. Hence, all the regions $\mathcal{T}(2^{n},1,k)$ are non-empty for $k\geqslant 4n+2$. Lemma is proved.

\begin{figure}[h]
\begin{minipage}[h]{1.0\linewidth}
\center{\includegraphics[width=0.7\linewidth]{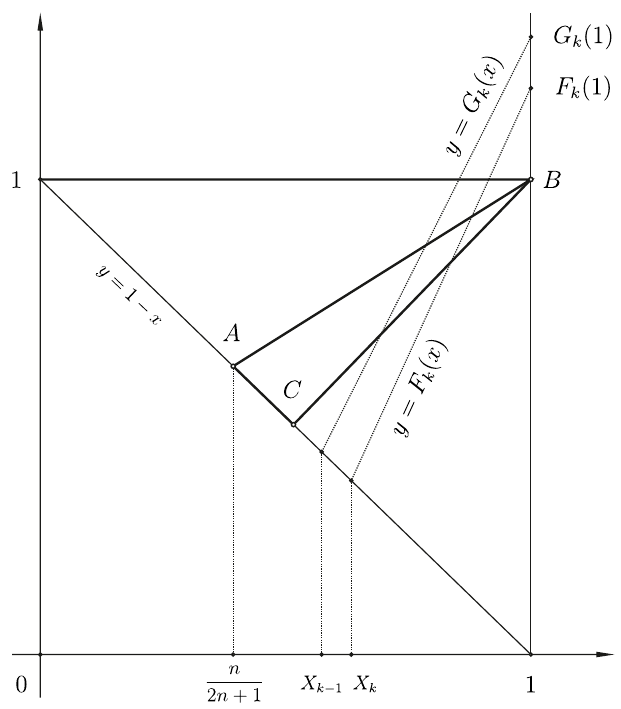}}\\
\end{minipage}
\begin{center}
Fig.~4b. \emph{Schematic drawing of the polygon $\mathcal{T}(2^{n},1)$ with the lines $y = F_{k}(x)$, $y = G_{k}(x)$}.
\end{center}
\end{figure}
\vspace{0.2cm}

\textsc{Lemma 11.} \emph{The region $\mathcal{T}(3,2,1,k)$ is non-empty iff  $7\leqslant k\leqslant 12$. The area of the union of non-empty regions $\mathcal{T}(3,2,1,k)$ is equal to} $1/35$.
\vspace{0.3cm}

\textsc{Proof.} The direct calculation shows that $\mathcal{T}(3,2,1)$ is the quadrangle with the vertices
\[
A = \biggl(\frac{4}{7},\frac{3}{7}\biggr),\quad B = \biggl(1,\frac{3}{5}\biggr), \quad C = \biggl(1,\frac{4}{7}\biggr), \quad D = \biggl(\frac{3}{5},\frac{2}{5}\biggr).
\]
It is described by the inequalities
\begin{equation}\label{lab_16}
\begin{cases}
\displaystyle \frac{4}{7}<x\leqslant \frac{3}{5},\quad 1-x<y\leqslant \frac{1+2x}{5},\\[12pt]
\displaystyle \frac{3}{5}<x\leqslant 1,\quad \frac{1+3x}{7}<y\leqslant \frac{1+2x}{5}
\end{cases}
\end{equation}
(see fig. 5) and its area is equal to $1/35$. When passing to the region $\mathcal{T}(3,2,1,k)$, the inequalities $f(x)<y\leqslant g(x)$,
\begin{align*}
& f(x) = f(x;3,2,1,k) = \frac{1+x\,\mathbb{K}_{3}(2,1,k+1)}{\mathbb{K}_{4}(3,2,1,k+1)} = \frac{1+x(k-1)}{2k-3},\\[12pt]
& g(x) = f(x;3,2,1,k) = \frac{1+x\,\mathbb{K}_{3}(2,1,k)}{\mathbb{K}_{4}(3,2,1,k)} = \frac{1+x(k-2)}{2k-5},
\end{align*}
are added to the system (\ref{lab_16}). The solutions of the inequality $g(x)\leqslant (1+3x)/7$ form the ray $x\leqslant (2k-12)/(k+1)$/ For $k\geqslant 13$, it contains the segment $0\leqslant x\leqslant 1$. Hence, if $k\geqslant 13$ then the quadrangle $\mathcal{T}(3,2,1)$ lies over the line $y = g(x)$. Thus, the region $\mathcal{T}(3,2,1,k)$ is empty for $k\geqslant 13$.

Further, the solutions of the inequality $f(x)\geqslant (1+2x)/5$ form the ray
$$
x\geqslant \frac{2k-8}{k+1}.
$$
If $k\leqslant 6$, then it contains the segment $4/7\leqslant x\leqslant 1$. Hence, for $k\leqslant 6$, the quadrangle $\mathcal{T}(3,2,1)$ lies entirely under the line $y = f(x)$. This implies the emptiness of the regions $\mathcal{T}(3,2,1,k)$, $1\leqslant k\leqslant 6$.

The non-emptiness of the regions $\mathcal{T}(3,2,1,k)$, $7\leqslant k\leqslant 12$, is checked directly (see fig.~5).

\vspace{0.3cm}

\textsc{Remark.} Both the precise form of the inequalities describing all the regions $\mathcal{T}(3,2,1,k)$, $7\leqslant k\leqslant 12$, and the values of its areas are given in Appendix I.
\vspace{0.3cm}

\textsc{Lemma 12.} \emph{Suppose that $n\geqslant 2$. Then the region $\mathcal{T}(3,2^{n},1,k)$ is non-empty only for $4n+2\leqslant k\leqslant 4n+8$.
The union of all non-empty regions $\mathcal{T}(3,2^{n},1,k)$ has an area equal to}
\[
\frac{2}{(2n+1)(2n+3)(2n+5)}.
\]
\vspace{0.3cm}

\begin{figure}[h]
\begin{minipage}[h]{1.0\linewidth}
\center{\includegraphics[width=0.65\linewidth]{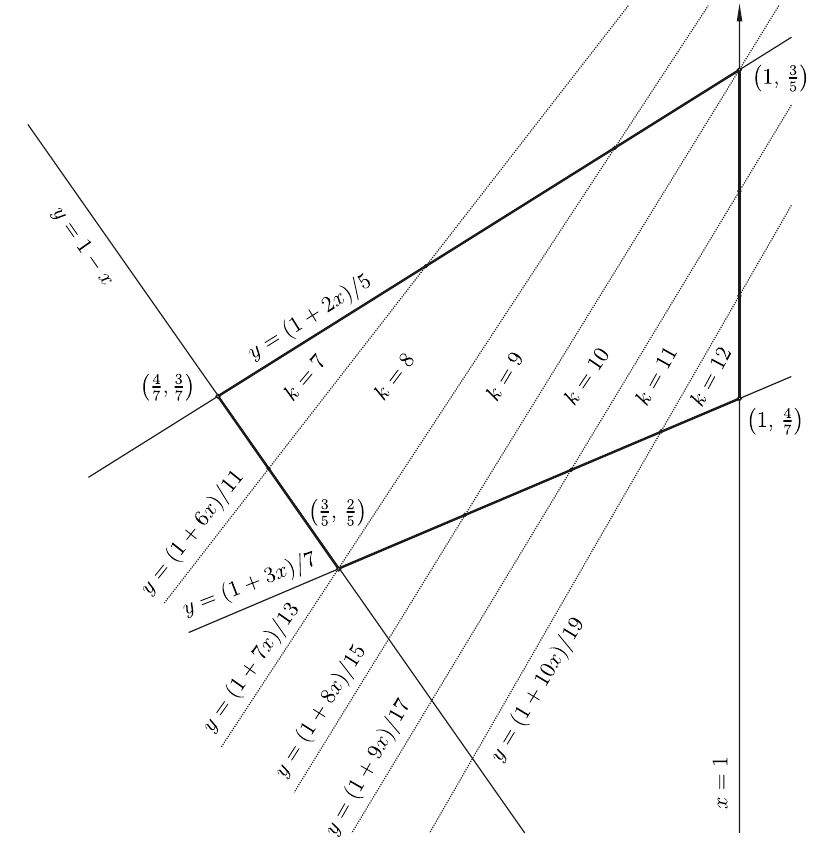}}\\
\end{minipage}
\begin{center}
Fig.~5. \emph{Schematic drawing of the polygon $\mathcal{T}(3,2,1)$}.
\end{center}
\end{figure}

\textsc{Proof.} Since $\mathbb{K}_{n+3}(3,2^{n},1,k) = 2k-2n-3<0$ for $k\leqslant n+1$, then the region $\mathcal{T}(3,2^{n},1,k)$ is empty for such $k$. In what follows, we assume that $k\geqslant n+2$.

Lemma 9(a) together with the system (\ref{lab_14}) imply that the region $\mathcal{T}(3,2^{n},1)$ is described by the system
\begin{equation}\label{lab_17}
\begin{cases}
\displaystyle \frac{2n-1}{2n+1} < x\leqslant \frac{2n+1}{2n+3},\quad \frac{1+x}{4} < y\leqslant \frac{1+(n+1)x}{2n+3},\\[12pt]
\displaystyle \frac{2n+1}{2n+3} < x\leqslant 1,\quad \frac{1+(n+2)x}{2n+5} < y\leqslant \frac{1+(n+1)x}{2n+3}
\end{cases}
\end{equation}
and is the quadrangle with the vertices
\begin{multline*}
A = \biggl(\frac{2n-1}{2n+1},\frac{n}{2n+1}\biggr),\quad B = \biggl(1,\frac{n+2}{2n+3}\biggr),\quad C = \biggl(1,\frac{n+3}{2n+5}\biggr),\\
D = \biggl(\frac{2n+1}{2n+3},\frac{n+1}{2n+3}\biggr).
\end{multline*}
When passing to the region $\mathcal{T}(3,2^{n},1,k)$, the inequalities $F_{k}(x)<y\leqslant G_{k}(x)$,
\begin{align*}
& F_{k}(x) = f_{n+3}(x;3,2^{n},1,k) = \frac{1+x\,\mathbb{K}_{n+2}(2^{n},1,k+1)}{\mathbb{K}_{n+3}(3,2^{n},1,k+1)} = \frac{1+x(k-n)}{2(k-n)-1}, \\[12pt]
& G_{k}(x) = g_{n+3}(x;3,2^{n},1,k) = \frac{1+x\,\mathbb{K}_{n+2}(2^{n},1,k)}{\mathbb{K}_{n+3}(3,2^{n},1,k)} = \frac{1+x(k-n-1)}{2(k-n)+1}.
\end{align*}
are added to the system (\ref{lab_17}).
The solutions of the inequality
\[
\frac{1+(n+1)x}{2n+3}\leqslant F_{k}(x)
\]
form the ray $x\geqslant x_{k}$, where
\[
x_{k} = \frac{2(k-2n-2)}{k+1}.
\]
If $k\leqslant 4n+1$ then it is easy to check that $x_{k}$ does not exceed the $x$-coordinate of the vertex $A$, that is,
\[
\frac{2(k-2n-2)}{k+1} \leqslant \frac{2n-1}{2n+1}.
\]
Hence, for such $k$, the ray $x\geqslant x_{k}$ contains the segment
\begin{equation}\label{lab_18}
\frac{2n-1}{2n+1}\leqslant x\leqslant 1,
\end{equation}
which is the projection of $ABCD$ onto horizontal axis. Therefore, the quadrangle $ABCD$ lies entirely under the line $y = F_{k}(x)$ and the region
$\mathcal{T}(3,2^{n},1,k)$ is empty.

One can check similarly, that the solutions of the inequality
\[
G_{k}(x)\leqslant \frac{1+(n+2)x}{2n+5}
\]
form the ray $x\leqslant X_{k}$, where
\[
X_{k} = \frac{2(k-2n+4)}{k+1}.
\]
However, $X_{k}\geqslant 1$ for any $k\geqslant 4n+9$. Hence, such ray contains the segment (\ref{lab_18}). Therefore, $ABCD$ lies entirely over the line $y = G_{k}(x)$ and the region
$\mathcal{T}(3,2^{n},1,k)$ is empty.

The non-emptiness of the regions $\mathcal{T}(3,2^{n},1,k)$, $4n+2\leqslant k\leqslant 4n+8$ is checked directly (see fig. 6). Lemma is proved.
\vspace{0.3cm}

\textsc{Remark.} Both the precise form of the inequalities describing all the regions $\mathcal{T}(3,2^{n},1,k)$, $4n+2\leqslant k\leqslant 4n+8$, and the values of its areas are given in Appendix II.
\vspace{0.3cm}

\begin{figure}[h]
\begin{minipage}[h]{1.0\linewidth}
\center{\includegraphics[width=0.95\linewidth]{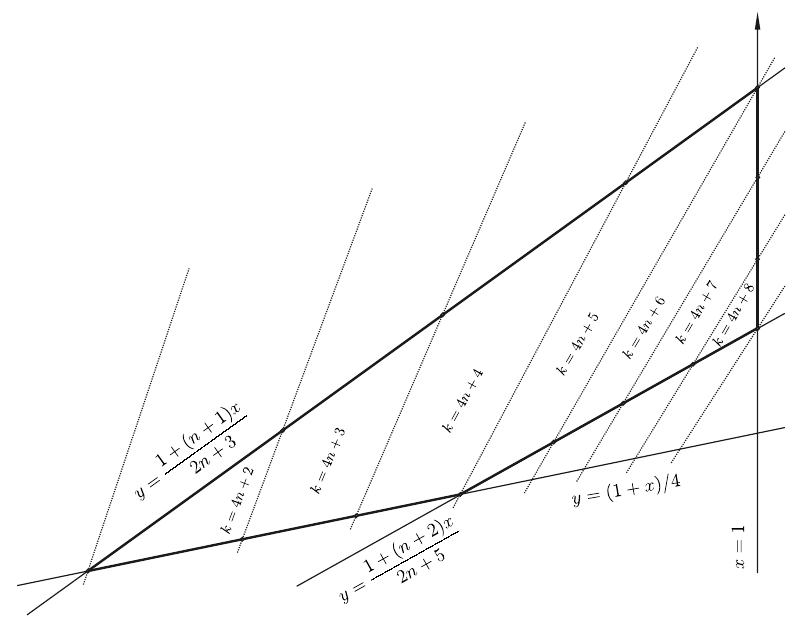}}\\
\end{minipage}
\begin{center}
Fig.~6. \emph{Schematic drawing of the region $\mathcal{T}(3,2^{n},1)$ and its subregions $\mathcal{T}(3,2^{n},1,k)$, $4n+2\leqslant k\leqslant 4n+8$, for $n\geqslant 2$}.
\end{center}
\end{figure}

\section{Iterative construction of the sets $\boldsymbol{\mathcal{A}_{r}^{\circ}(3,1,2)}$}
\vspace{0.5cm}

Together with the sets $\mathcal{A}_{r}^{\circ}(D,c_{0},c_{1})$, it is convenient to introduce the sets $\mathcal{A}_{r}^{*}(D,c_{0},c_{1})$.
Their definition differs only in that the last congruence in (\ref{lab_07}) is replaced by its negation:
\[
c_{1}\,\mathbb{K}_{r}(k_{1},\ldots,k_{r})-c_{0}\,\mathbb{K}_{r-1}(k_{2},\ldots,k_{r})\not\equiv c_{0}\pmod{D}.
\]
Hence, the condition for the tuple to belong to the set $\mathcal{A}_{r}^{*}(3,1,2)$ can be written in the form
\begin{equation}\label{lab_19}
\mathbb{K}_{i}(k_{1},\ldots,k_{i})+\mathbb{K}_{i-1}(k_{2},\ldots,k_{i})\not\equiv 2\pmod{3},\quad i = 1,2,\ldots, r.
\end{equation}

\textsc{Case} $r=1$. Condition (\ref{lab_09}) takes the form
\[
k_{1}+1\equiv 2\pmod{3},\quad\text{that is,}\quad k_{1}\equiv 1\pmod{3}.
\]
Thus the set $\mathcal{A}_{1}^{\circ}(3,1,2)$ consists of all tuples $(k_{1})$, where $k_{1}\geqslant 1$, $k_{1}\equiv 1\pmod{3}$, and the set $\mathcal{A}_{1}^{*}(3,1,2)$ is formed by al tuples $(k_{1})$, $k_{1}\geqslant 2$, $k_{1}\not\equiv 1\pmod{3}$, correspondingly.
\vspace{0.3cm}

\textsc{Case} $r=2$. Condition (\ref{lab_09}) takes the form
\begin{equation*}
\begin{cases}
k_{1}\not\equiv 1\pmod{3},\\
k_{2}(k_{1}+1)\equiv 0\pmod{3}.
\end{cases}
\end{equation*}
If $k_{1}\equiv 0\pmod{3}$, then we necessary have $k_{2}\equiv 0\pmod{3}$ and therefore $k_{1}\geqslant 3$, $k_{2}\geqslant 3$. By lemma 2(c), such regions are empty. Suppose that $k_{1}\equiv 2\pmod{3}$. Then $k_{2}>0$ is arbitrary integer. Now let us peak out the values corresponding to non-empty regions $\mathcal{T}(k_{1},k_{2})$.

If $k_{1}\ne 2$ then $k_{1}\geqslant 5$. By lemma 2(d) and lemma 4, in this case, there is a unique non-empty region $\mathcal{T}(k_{1},1)$. If $k_{1} = 2$ then, by lemma 2(d), the regions $\mathcal{T}(k_{1},k_{2}) = \mathcal{T}(2,k_{2})$ are non-empty only for $1\leqslant k_{2}\leqslant 4$. Thus, the set $\mathcal{A}_{2}^{\circ}(3,1,2)$ consists of the tuples
\begin{align}\label{lab_20}
& (k_{1},1),\quad\text{where}\quad k_{1}\equiv 2\pmod{3},\quad k_{1}\geqslant 2;\\
& (2,2), (2,3), (2,4).\notag
\end{align}
Further, the conditions (\ref{lab_19}) take the form
\begin{equation*}
\begin{cases}
k_{1}\not\equiv 1\pmod{3},\\
k_{2}(k_{1}+1)\not\equiv 0\pmod{3}
\end{cases}
\quad\text{or, that is the same,}\quad
\begin{cases}
k_{1} \equiv 0\pmod{3},\\
k_{2}\not\equiv 0\pmod{3}.
\end{cases}
\end{equation*}
By lemma 2(c),(d) and lemma 4, we conclude that, in the case $k_{1} = 3$, such region is non-empty only for $k_{2}=1,2$, and, in the case $k_{1}\equiv 0\pmod{3}$, $k_{1}\geqslant 6$, is non-empty only for $k_{1}=1$. Thus, $\mathcal{A}_{2}^{*}(3,1,2)$ consists of the tuples
\begin{align}\label{lab_21}
& (k_{1},1),\quad\text{where}\quad k_{1}\equiv 0\pmod{3},\quad k_{1}\geqslant 3;\\
& (3,2).
\end{align}

\textsc{Case} $r=3$. Any triple $(k_{1},k_{2},k_{3})\in \mathcal{A}_{3}^{\circ}(3,1,2)$ necessarily satisfies to the condition $(k_{1},k_{2})\in \mathcal{A}_{2}^{*}(3,1,2)$. By (\ref{lab_21}),
we either have $k_{1}\equiv 0\pmod{3}$, $k_{1}\geqslant 3$, $k_{2}=1$, or $k_{1} = 3, k_{2} = 2$.

In the first case, the condition (\ref{lab_09}) corresponding to $i=r=3$ takes the form
\[
(k_{1}k_{3}-k_{1}-k_{3})+(k_{1}-1)\equiv 2\pmod{3},\quad\text{т.е.}\quad k_{1}(k_{3}-1)\equiv 0\pmod{3}.
\]
One can easily see that it is satisfied automatically for any  $k_{3}\geqslant 1$. Using lemma 3 (2b, 2e, 2f), we conclude that non-empty regions correspond to the following triples:
\begin{align}\label{lab_22}
& (k_{1},1,2),\quad k_{1}\equiv 0\pmod{3},\quad k_{1}\geqslant 6;\notag \\
& (3,1,4),\; (3,1,5),\; (3,1,5),\; (3,1,7),\; (3,1,8);\\
& (6,1,3).\notag
\end{align}
In the second case, for the triple $(k_{1},k_{2},k_{3}) = (3,2,k_{3})$, the condition (\ref{lab_09}) has the form
\[
(5k_{3}-3)+(2k_{3}-1)\equiv 2\pmod{3},\quad\text{that is,}\quad 7k_{3}\equiv 0\pmod{3}.
\]
But the emptiness of the regions $\mathcal{T}(3,2,k_{3})$ for $k_{3}\geqslant 3$ follows from lemma 3(3c).
Therefore,, $\mathcal{A}_{3}^{\circ}(3,1,2)$ consists of the triples (\ref{lab_22}).

Further, if $(k_{1},k_{2},k_{3})\in \mathcal{A}_{3}^{*}(3,1,2)$, then we necessarily have $(k_{1},k_{2})\in \mathcal{A}_{2}^{*}(3,1,2)$. As it was stated before, we either have $k_{1}\equiv 0\pmod{3}, k_{1}\geqslant 3$, $k_{2}=1$, or $k_{1}=3$, $k_{2}=2$.

For $k_{2}=1$, the condition (\ref{lab_19}) corresponding to $i=r=3$, reduces to
\[
k_{1}(k_{3}-1)\not\equiv 0\pmod{3}.
\]
Since $k_{1}\equiv 0\pmod{3}$ in the first case, this condition can not be satisfied for any $k_{3}$. In the second case, the same condition takes the form
\[
7k_{3}\not\equiv 0\pmod{3}.
\]
For $k_{3}\geqslant 3$, the emptiness of $\mathcal{T}(3,2,k_{3})$ follows from lemma 3 (3c). Therefore, non-empty regions correspond to the values  $k_{3}=1, k_{3} = 2$ only. Hence, the set $\mathcal{A}_{3}^{*}(3,1,2)$ consists of the triples
\begin{equation}\label{lab_23}
(3,2,1),\quad (3,2,2).
\end{equation}

It turns out that for any $r\geqslant 4$ the set $\mathcal{A}_{r}^{*}(3,1,2)$ consists of two tuples: $(3,2^{r-2},1)$ and $(3,2^{r-1})$.
Let us prove it by induction with respect to $r$.

For $r=3$, it was proved earlier. Suppose now that this assertion holds for any $r$, $3\leqslant r\leqslant m-1$, and check it for $r=m$.

Indeed, let $(k_{1},\ldots,k_{m})\in \mathcal{A}_{m}^{*}(3,1,2)$. Then $(k_{1},\ldots,k_{m-1})\in \mathcal{A}_{m-1}^{*}(3,1,2)$, so, by the induction, $(k_{1},\ldots,k_{m-1})$ is either $(3,2^{m-3},1)$ or $(3,2^{m-2})$. Hence, $(k_{1},\ldots,k_{m})$ has the form  $(3,2^{m-3},1,k)$ or $(3,2^{m-2},k)$, where $k\geqslant 1$ is some integer.

In the first case, the condition (\ref{lab_19}) corresponding to $i=m$, has the form
\begin{equation}\label{lab_24}
\mathbb{K}_{m}(3,2^{m-3},1,k)+\mathbb{K}_{m-1}(2^{m-3},1,k)\not\equiv 2\pmod{3}.
\end{equation}
By explicit formulas for the continuants, we get
\[
\mathbb{K}_{m}(3,2^{m-3},1,k) = 2(k-m)+3,\quad \mathbb{K}_{m-1}(2^{m-3},1,k) = k-m+2,
\]
so (\ref{lab_24}) has the form $3(k-m)+5\not\equiv 2\pmod{3}$. Therefore, it is not satisfied for any $k$. Hence, $(k_{1},\ldots,k_{m})$ is not of the form $(3,2^{m-3},1,k)$.

In the second case, we have the condition
\[
\mathbb{K}_{m}(3,2^{m-2},k) + \mathbb{K}_{m-1}(2^{m-2},k) \not\equiv 2\pmod{3}.
\]
Using the formulas
\[
\mathbb{K}_{m}(3,2^{m-2},k) = 2km-k-2m+3,\quad \mathbb{K}_{m}(2^{m-2},k) = km-k-m+2,
\]
we rewrite it as follows
\[
3m(k-1)-2k+5\not\equiv 2\pmod{3},\quad\text{that is,}\quad k\not\equiv 0\pmod{3}.
\]
Since $m-2\geqslant 2$ then lemma 9(b) implies that the region $\mathcal{T}(3,2^{m-2},k)$ is empty for $k\geqslant 3$ and is non-empty for $k = 1,2$. Thus, $\mathcal{A}_{m-1}^{*}(3,1,2)$ consists of two tuples $(3,2^{m-2},1)$ and $(3,2^{m-1})$.

Now let us find all the tuples that contain the set $\mathcal{A}_{r}^{\circ}(3,1,2)$ for $r\geqslant 4$.

If $(k_{1},\ldots,k_{r})$ lies in $\mathcal{A}_{r}^{\circ}(3,1,2)$, then $(k_{1},\ldots,k_{r-1})$ belongs to $\mathcal{A}_{r-1}^{*}(3,1,2)$ and hence coincides either with с $(3,2^{r-2})$ or with $(3,2^{r-3},1)$.

In the first case, the initial tuple has the form $(3,2^{r-2},k)$, where $k$ satisfies to the last of conditions (\ref{lab_09}), that is,
\[
\mathbb{K}_{r}(3,2^{r-2},k) + \mathbb{K}_{r-1}(2^{r-2},k)\equiv 2 \pmod{3}
\]
or, that is the same, to the congruence
\[
3(k-1)(r-2)+4k-1\equiv 2\pmod{3},
\]
which is valid for $k\equiv 0\pmod{3}$. The emptiness of such region follows from lemma 9 (3b).

In the second case, the initial tuple has the form $(3,2^{r-2},1,k)$, so the condition for $k$ looks as
\[
\mathbb{K}_{r}(3,2^{r-3},1,k) + \mathbb{K}_{r-1}(2^{r-3},1,k)\equiv 2 \pmod{3}.
\]
The last congruence is equivalent to the following:
\[
(2k-2r+3)+(k-r+2)\equiv 2 \pmod{3}.
\]
However, it holds for any $k\geqslant 1$. If $r = 4$, then, by lemma 11, non-empty regions $\mathcal{T}(3,2^{r-3},1,k) = \mathcal{T}(3,2,1,k)$ corresponds only to
 $7\leqslant k\leqslant 12$. In the case $r\geqslant 5$, lemma 12 implies that such non-empty regions correspond only to the values $4r-10\leqslant k\leqslant 4r-4$.

This completes the description of the sets  $\mathcal{A}_{r}^{\circ}(3,1,2)$.

It is convenient to formulate the above statements as follows.
\vspace{0.3cm}

\textsc{Lemma 13.} \emph{The set} $\mathcal{A}_{r}^{\circ}(3,1,2)$ (\emph{and therefore the set} $\mathcal{A}_{r}^{\circ}(3,2,1)$) \emph{consists of the following tuples}:

\begin{longtable}{|>{\fontsize{10}{9pt}\selectfont}c|>{\fontsize{10}{9pt}\selectfont}c|>{\fontsize{10}{9pt}\selectfont}c|}
\hline $r$ & $\mathbf{k} = (k_{1},\ldots,k_{r})$       & $\mathbb{K}_{r}(\mathbf{k})$  \\
\hline 1   & $(k_{1}), k_{1}\equiv 1\pmod{3}$, $k_{1}\geqslant 3$ & $k_{1}$\\
\hline 2   & $(k_{1},1), k_{1}\equiv 2\pmod{3}$, $k_{1}\geqslant 2$ & $k_{1}-1$\\
\hline 2   & $(2,2)$ & $3$\\
\hline 2   & $(2,3)$ & $5$\\
\hline 2   & $(2,4)$ & $7$\\
\hline 3   & $(k_{1},1,2), k_{1}\equiv 0\pmod{3}$, $k_{1}\geqslant 3$ & $k_{1}-2$\\
\hline 3   & $(3,1,k_{3})$, $4\leqslant k_{3}\leqslant 8$ & $2k_{3}-3$\\
\hline 3   & $(6,1,3)$ & $9$\\
\hline 4   & $(3,2,1,k_{4})$, $7\leqslant k_{4}\leqslant 12$ & $2k_{4}-5$\\
\hline $\geqslant 5$   & $\bigl(3,2^{r-3},1,k_{r}\bigr)$, $4r-10\leqslant k_{r}\leqslant 4r-4$ & $2k_{r}-2r+3$\\
\hline
\end{longtable}
\vspace{0.5cm}

\section{Iterative construction of the sets $\boldsymbol{\mathcal{A}_{r}^{\circ}(3,1,0)}$}
\vspace{0.5cm}

The conditions for the tuple $(k_{1},\ldots,k_{r})$ to belong to the set $\mathcal{A}_{r}^{\circ}(3,1,0)$ has the form (\ref{lab_10}).
Hence, the condition for the tuple to belong to the set $\mathcal{A}_{r}^{*}(3,1,0)$ is the following:
\begin{equation}\label{lab_25}
\mathbb{K}_{i-1}(k_{2},\ldots,k_{i})\not\equiv 2\pmod{3},\quad i = 2,3,\ldots, r-1.
\end{equation}
As was mentioned earlier, in the \textsc{case}  $r=1$, the conditions (\ref{lab_10}) can not be satisfied for any $k_{1}$. So, the set $\mathcal{A}_{1}^{\circ}(3,1,0)$ is empty, and the set
$\mathcal{A}_{1}^{*}(3,1,0)$ consists of all unit tuples $(k_{1})$, $k_{1}\geqslant 1$.

In the \textsc{case} $r=2$, the condition (\ref{lab_10}) takes the form $k_{2}\equiv 2\pmod{3}$. Among the tuples $(k_{1},k_{2})$ with $k_{1}\geqslant 1$, $k_{2}\equiv 2\pmod{3}$, only the pairs
\begin{align*}
& (1,k_{2}),\quad k_{2}\equiv 2\pmod{3},\quad k_{2}\geqslant 2;\\
& (2,2),\; (3,2),\; (4,2),
\end{align*}
correspond to non-empty regions. All these pairs form the set $\mathcal{A}_{2}^{\circ}(3,1,0)$.
Further, the conditions (\ref{lab_25}) take the form $k_{2}\not\equiv 2\pmod{3}$, $k_{1}\geqslant 1$, where $k_{1}\geqslant 1$ is arbitrary integer.
Using lemma 2 (b),(d), in the cases $k_{1}=1$ and $k_{2}=1$, we get the following pairs corresponding to non-empty regions:
\begin{align*}
& (1,k_{2}),\quad k_{2}\not\equiv 2\pmod{3},\quad k_{2}\geqslant 3\\
& (k_{1},1),\quad k_{1}\geqslant 2\quad - \quad\text{arbitrary}.
\end{align*}

Due to the above arguments, in the case $k_{1} = 2$ we have only two non-empty regions: $\mathcal{T}(2,3)$ and $\mathcal{T}(2,4)$.
In the case $k_{1}\geqslant 3$, there are no non-empty regions corresponding to the values $k_{2}\geqslant 3$, $k_{2}\not\equiv 2\pmod{3}$.
Hence, $\mathcal{A}_{2}^{*}(3,1,0)$ is contained of the pairs
\begin{align}\label{lab_26}
& \text{(a)}\quad (1,k_{2}),\quad k_{2}\equiv 0\pmod{3},\quad k_{2}\geqslant 3;\notag \\
& \text{(b)}\quad (1,k_{2}),\quad k_{2}\equiv 1\pmod{3},\quad k_{2}\geqslant 4; \\
& \text{(c)}\quad (k_{1},1),\quad k_{1}\geqslant 2; \notag \\
& \text{(d)}\quad (2,3),\;(2,4).\notag
\end{align}
In the \textsc{case} $r=3$, the set $\mathcal{A}_{r}^{\circ}(3,1,0)$ consists of the triples $(k_{1},k_{2},k_{3})$, where $(k_{1},k_{2})\in \mathcal{A}_{2}^{*}(3,1,0)$ is one of the pairs (\ref{lab_26}) and $k_{3}$ satisfies to the condition
\begin{equation}\label{lab_27}
k_{2}k_{3}\equiv 0\pmod{3}.
\end{equation}
If $(k_{1},k_{2})$ has the form (\ref{lab_26}\,a), then (\ref{lab_27}) holds for any $k_{3}$. By lemma 3 (4a), (6), only the triples
\begin{align*}
& (1,k_{2},1),\quad k_{2}\equiv 0\pmod{3},\quad k_{2}\geqslant 3;\\
& (1,3,2)
\end{align*}
correspond to non-empty regions.
If $(k_{1},k_{2})$ has the form (\ref{lab_26}\,b), then, by (\ref{lab_27}), we have $k_{3}\equiv 0\pmod{3}$. Since $k_{2}\geqslant 4$, then, by lemma 2 (c), all the corresponding regions are empty.

If $(k_{1},k_{2})$ has the form (\ref{lab_26}\,c), then we similarly have $k_{3}\equiv 0\pmod{3}$. By lemma 3, (2a)-(2f), only the triples
\begin{align*}
& (2,1,k_{3}),\quad k_{3}\equiv 0\pmod{3},\quad k_{3}\geqslant 6;\\
& (3,1,6);\\
& (k_{1},1,3),\quad 4\leqslant k_{1}\leqslant 8
\end{align*}
correspond to non-empty regions.

If $(k_{1},k_{2}) = (2,3)$, then, in view of (\ref{lab_27}), $k_{3}$ can be arbitrary. However, by lemma 3 (4b), non-empty regions correspond only to triples
\[
(2,3,1),\quad (2,3,2).
\]
Finally, if $(k_{1}, k_{2}) = (2,4)$ then we necessarily have $k_{3}\equiv 0\pmod{3}$ and hence $k_{3}\geqslant 3$. By lemma 2, all the corresponding regions are empty.

Thus, $\mathcal{A}_{3}^{\circ}(3,1,0)$ consists of the triples
\begin{align}\label{lab_28}
& (1,k_{2},1),\quad k_{2}\equiv 0\pmod{3},\quad k_{2}\geqslant 3;\notag \\
& (2,1,k_{3}),\quad k_{3}\equiv 0\pmod{3},\quad k_{3}\geqslant 6;\notag\\
& (1,3,2);\\
& (2,3,1),\;(2,3,2);\notag\\
& (3,1,6);\notag\\
& (k_{1},1,3),\quad 4\leqslant k_{1}\leqslant 8.\notag
\end{align}

Passing to the description of the triples $(k_{1},k_{2},k_{3})\in \mathcal{A}_{3}^{*}(3,1,0)$, in view of (\ref{lab_25}), (\ref{lab_27}), one has
\begin{equation}\label{lab_29}
k_{2}k_{3}\not\equiv 0\pmod{3}.
\end{equation}
The condition excludes from consideration both the pairs (\ref{lab_26}\,a) and the pair $(2,3)$.

Let $(k_{1},k_{2})$ be of the form (\ref{lab_26}\,b). Here $k_{2}\geqslant 4$, so, in view of lemma 2 (c),(d) and lemma 4, it is sufficient to consider only the values $k_{3} = 1,2$. By lemma 3 (5a),(6), only the triples
\begin{align*}
& (1,k_{2},1),\quad k_{2}\equiv 1\pmod{3},\quad k_{2}\geqslant 4;\\
& (1,4,2)
\end{align*}
correspond to non-empty regions.
Suppose that $(k_{1},k_{2})$ has the form (\ref{lab_26}\,c). Then $k_{3}\not\equiv 0\pmod{3}$ and hence, by lemma 3 (2a)-(2f), non-empty regions correspond only to the triples
\begin{align*}
& (2,1,k_{3}),\quad k_{3}\equiv 1\pmod{3},\quad k_{3}\geqslant 7;\\
& (2,1,k_{3}),\quad k_{3}\equiv 2\pmod{3},\quad k_{3}\geqslant 8;\\
& (3,1,k_{3}),\quad k_{3}=4,5,7,8;\\
& (4,1,4),\;(4,1,5);\\
& (5,1,4);\\
& (k_{1},1,2),\quad k_{1}\geqslant 6.
\end{align*}
Finally, let $(k_{1},k_{2}) = (2,4)$. Then the single non-empty region corresponds to the triple $(2,4,1)$.

Combining the obtained results, we conclude that $\mathcal{A}_{3}^{*}(3,1,0)$ consists of the following triples:
\begin{align}\label{lab_30}
& (k_{1},1,2),\quad k_{1}\geqslant 6;\notag \\
& (1,k_{2},1),\quad k_{2}\equiv 1\pmod{3},\quad k_{2}\geqslant 4;\notag\\
& (2,1,k_{3}),\quad k_{3}\equiv 1\pmod{3},\quad k_{3}\geqslant 7;\notag\\
& (2,1,k_{3}),\quad k_{3}\equiv 2\pmod{3},\quad k_{3}\geqslant 8;\notag\\
& (1,4,2),\;(2,4,1);\\
& (3,1,k_{3}),\quad k_{3} = 4,5,7,8;\notag\\
& (4,1,4),\;(4,1,5);\notag \\
& (5,1,4).\notag
\end{align}

In the \textsc{case} $r=4$, for any quadruple $(k_{1},k_{2},k_{3},k_{4})\in \mathcal{A}_{4}^{\circ}(3,1,0)$, one has: $(k_{1},k_{2},k_{3})\in \mathcal{A}_{3}^{*}(3,1,0)$ is the triple from the list (\ref{lab_30}) and the value $k_{4}$ satisfies to the following condition:
\begin{equation}\label{lab_31}
k_{4}(k_{2}k_{3}-1) - k_{2}\equiv 2\pmod{3}.
\end{equation}
If $(k_{1},k_{2},k_{3}) = (k_{1},1,2)$, then $k_{4}\equiv 0\pmod{3}$. By lemma 3 (3a), the unique value of $k_{4}$ corresponding to non-empty region, is $3$. By lemma 4, both the regions $\mathcal{T}(k,1,2,3)$ and $\mathcal{T}(3,2,1,k)$ are simultaneously empty or non-empty. By lemma 11, the region $\mathcal{T}(3,2,1,k)$ is non-empty only for $7\leqslant k\leqslant 12$. Hence, only the following quadruples belong to the set $\mathcal{A}_{4}^{\circ}(3,1,0)$:
\[
(k_{1},1,2,3),\quad 7\leqslant k_{1}\leqslant 12.
\]

If $(k_{1},k_{2},k_{3}) = (1,k_{2},1)$, where $k_{2}\equiv 1\pmod{3}$, $k_{2}\geqslant 4$, then (\ref{lab_31}) holds for any $k_{4}$.
By lemma 3 (2c), (2e), (2f), the necessary condition, that is, the non-emptiness of $\mathcal{T}(k_{2},1,k_{4})$, implies $3\leqslant k_{4}\leqslant 5$ for $k_{2} = 4$, implies
$2\leqslant k_{4}\leqslant 3$ for $k_{2} = 7$ and implies $k_{4} = 2$ for $k_{2}\geqslant 11$.

In first two cases, the direct calculation shows that non-empty regions correspond to the quadruples
\begin{align*}
& (1,4,1,3),\;(1,4,1,4),\;(1,4,1,5);\\
& (1,7,1,3);\\
& (1,k_{2},1,2),\quad k_{2}\equiv 1\pmod{3},\quad k_{2}\geqslant 7.
\end{align*}
In the third case, the non-emptiness follows from lemma 7 (d) for $k_{2}\geqslant 9$ and is checked directly for $k_{2} = 7,8$.

If $(k_{1},k_{2},k_{3}) = (2,1,k_{3})$, where $k_{3}\equiv 1\pmod{3}$, $k_{3}\geqslant 7$, then (\ref{lab_31}) also holds for any $k_{4}\geqslant 1$. By lemma 3 (6), the
necessary condition, that is, non-emptiness of the region $\mathcal{T}(k_{2},k_{3},k_{4}) = \mathcal{T}(1,k_{3},k_{4})$, implies the equality $k_{4}=1$.
By lemma 7 (d), all the regions $\mathcal{T}(2,1,k_{3},1)$ are non-empty for $k_{3}\geqslant 9$. For $k_{3} = 7,8$, their non-emptiness can be checked directly. Thus we get the quadruples
\[
(2,1,k_{3},1),\quad k_{3}\equiv 1\pmod{3},\quad k_{3}\geqslant 7.
\]
If $(k_{1},k_{2},k_{3}) = (2,1,k_{3})$, where $k_{3}\equiv 2\pmod{3}$, $k_{3}\geqslant 8$, then it follows from (\ref{lab_31}) that $k_{4}\equiv 0\pmod{3}$ and, in particular, $k_{4}\geqslant 3$. Since $k_{3}\geqslant 8$, then the necessary condition, that is, non-emptiness of  $\mathcal{T}(k_{3},k_{4})$, fails by lemma 2.  Therefore, all the regions under considering are empty.

By lemma 3 (5b), the region $\mathcal{T}(k_{2},k_{3},k_{4}) = \mathcal{T}(4,2,k_{4})$ is non-empty only for $k_{4}=1$. Hence, all the regions under considering are empty.

If $(k_{1},k_{2},k_{3}) = (2,4,1)$, then (\ref{lab_31}) is satisfied for an arbitrary $k_{4}\geqslant 1$. However, non-emptiness of the region $\mathcal{T}(4,1,k_{4})$
is the necessary condition for non-emptiness of the region $\mathcal{T}(k_{1},k_{2},k_{3},k_{4}) = \mathcal{T}(2,4,1,k_{4})$. By lemma 3 (2c), this implies the inequalities $3\leqslant k_{4}\leqslant 5$. The direct calculation shows that non-empty initial regions correspond to the quadruples $(2,4,1,3)$ and $(2,4,1,4)$.

If $(k_{1},k_{2},k_{3}) = (3,1,k_{3})$, where $k_{3}$ is one of the numbers $4,5,7,8$, then (\ref{lab_31}) takes the form
\begin{equation}\label{lab_32}
k_{4}(k_{3}-1)\equiv 0\pmod{3}.
\end{equation}
In the cases $k_{3} = 5,8$, this congruence implies the condition  $k_{4}\equiv 0\pmod{3}$ and, in particular, the inequality $k_{4}\geqslant 3$.
However, the non-emptiness of the region $\mathcal{T}(1,k_{3},k_{4})$ is necessary for the non-emptiness of $\mathcal{T}(k_{1},k_{2},k_{3},k_{4}) = \mathcal{T}(3,1,k_{3},k_{4})$.
By lemma 3 (5a),(6), such region is empty for $k_{4}\geqslant 2$. In the cases $k_{3} = 4, 7$, the congruence (\ref{lab_31}) is satisfied for any $k_{4}$.
However, lemma 3 (5a) implies that the region $\mathcal{T}(1,4,k_{4})$ is non-empty only for $k_{4} = 1,2$, and the region $\mathcal{T}(1,7,k_{4})$ is non empty only for $k_{4}=1$.
The direct calculation shows that all the corresponding regions $\mathcal{T}(3,1,k_{3},k_{4})$ are non-empty, and the corresponding quadruples are
\[
(3,1,4,1),\; (3,1,4,2),\; (3,1,7,1).
\]
If $(k_{1},k_{2},k_{3}) = (4,1,4)$, then (\ref{lab_31}) has the form $3k_{4}-1\equiv 2\pmod{3}$ and is satisfied for any $k_{4}\geqslant 1$.
The non-emptiness of the region $\mathcal{T}(1,4,k_{4})$ is a necessary condition for the non-emptiness of $\mathcal{T}(k_{1},k_{2},k_{3},k_{4}) = \mathcal{T}(4,1,4,k_{4})$.
Lemma 3 (5a) implies that $1\leqslant k_{4}\leqslant 2$. The direct calculation shows that each of there regions is non-empty.
The corresponding quadruples are:
\[
(4,1,4,1),\; (4,1,4,2).
\]
If $(k_{1},k_{2},k_{3}) = (4,1,5)$, then (\ref{lab_31}) has the form $4k_{4}-1\equiv 2\pmod{3}$ and is satisfied for any $k_{4}\equiv 0\pmod{3}$.
But the necessary condition for the non-emptiness of $\mathcal{T}(k_{1},k_{2},k_{3},k_{4}) = \mathcal{T}(4,1,5,k_{4})$, that is, the non-emptiness of the region
$\mathcal{T}(1,5,k_{4})$ fails by lemma 3 (6).

Finally, if $(k_{1},k_{2},k_{3}) = (5,1,4)$ then (\ref{lab_31}) has the form $3k_{4}-1\equiv 2\pmod{3}$ and is satisfied for any $k_{4}\geqslant 1$.
The non-emptiness of the region $\mathcal{T}(1,4,k_{4})$ is the necessary condition for the non-emptiness of $\mathcal{T}(k_{1},k_{2},k_{3},k_{4}) = \mathcal{T}(5,1,4,k_{4})$ and, by lemma 3 (5a), implies the inequalities $1\leqslant k_{4}\leqslant 2$. The direct calculation shows that the unique non-empty region corresponds to the quadruple
\[
(5,1,4,1).
\]
Therefore, the set $\mathcal{A}_{4}^{\circ}(3,1,0)$ consists of the quadruples
\begin{align}\label{lab_33}
& (1,k_{2},1,2),\quad k_{2}\equiv 1\pmod{3},\quad k_{2}\geqslant 7;\notag\\
& (2,1,k_{3},1),\quad k_{3}\equiv 1\pmod{3},\quad k_{3}\geqslant 7;\notag\\
& (1,4,1,3),\;(1,4,1,4),\;(1,4,1,5),\;(1,7,1,3);\notag\\
& (2,4,1,3),\;(2,4,1,4);\notag\\
& (3,1,4,1),\;(3,1,4,2),\;(3,1,7,1);\\
& (4,1,4,1),\;(4,1,4,2);\notag \\
& (5,1,4,1);\notag \\
& (k_{1},1,2,3),\quad 7\leqslant k_{1}\leqslant 12.\notag
\end{align}
Further, let us describe the set $\mathcal{A}_{4}^{*}(3,1,0)$. It consists of the quadruples $(k_{1},k_{2},k_{3},k_{4})$ such that $(k_{1},k_{2},k_{3})\in \mathcal{A}_{3}^{*}(3,1,0)$ and, moreover,
\begin{equation}\label{lab_34}
k_{4}(k_{2}k_{3}-1)-k_{2}\not\equiv 2\pmod{3}.
\end{equation}
First we note that if the congruence (\ref{lab_31}) is satisfied by any $k_{4}\geqslant 1$ then the congruence (\ref{lab_34}) fails for any $k_{4}$.
According to this, the cases when the triple $(k_{1},k_{2},k_{3})$ belongs to the list below, are immediately excluded from the consideration:
\begin{align*}
& (1,k_{2},1),\quad k_{2}\equiv 1\pmod{3},\quad k_{2}\geqslant 4;\\
& (2,1,k_{3}),\quad k_{3}\equiv 1\pmod{3},\quad k_{3}\geqslant 7;\\
& (2,4,1),\;(4,1,4),\;(5,1,4),\; (3,1,4),\; (3,1,7).
\end{align*}
In the rest cases, the condition (\ref{lab_31}) has the form $k_{4}\equiv 0\pmod{3}$, so the condition (\ref{lab_34}) is equivalent to the following:
\begin{equation}\label{lab_35}
k_{4}\not\equiv 0 \pmod{3}.
\end{equation}
If $(k_{1},k_{2},k_{3}) = (k_{1},1,2)$, where $k_{1}\geqslant 6$, then by lemma 3 (3a), the necessary condition for the non-emptiness of $\mathcal{T}(1,2,k_{4})$ leads to the inequalities
$2\leqslant k_{4}\leqslant 4$. Taking (\ref{lab_35}) into account, we obtain two possibilities: $k_{4} = 2$  $k_{4} = 4$.

In the first case, the non-empty regions correspond to the quadruples $(k_{1},1,2,2)$, $k_{1}\geqslant 10$, and in the second one - to the quadruples $(k_{1},1,2,4)$, $k_{1} = 6,7,8$.

If $(k_{1},k_{2},k_{3}) = (2,1,k_{3})$, $k_{3}\equiv 2\pmod{3}$, $k_{3}\geqslant 8$, then the necessary condition, that is, the non-emptiness of
the region $\mathcal{T}(k_{3},k_{4})$, gives $k_{4}=1$.
By lemma 7(d), for $k_{3}\geqslant 9$, all these regions are non-empty. In the case $k_{3} = 8$, the non-emptiness of $\mathcal{T}(2,1,k_{3},1)$ is checked directly. Therefore, all the quadruples
\[
(2,1,k_{3},1),\quad k_{3}\equiv 2\pmod{3},\quad k_{3}\geqslant 8,
\]
belong to $\mathcal{A}_{r}^{*}(3,1,0)$.

If $(k_{1},k_{2},k_{3}) = (1,4,2)$ then the necessary condition, that is, the non-emptiness of the region $\mathcal{T}(4,2,k_{4})$, leads by lemma 3 (3d) to the relation $k_{4}=1$ and, hence, to the single quadruple $(1,4,2,1)$.

If $(k_{1},k_{2},k_{3}) = (3,1,5)$ then the necessary condition, that is, the non-emptiness of $\mathcal{T}(1,5,k_{4})$ leaves by lemma 3 (6) the unique possibility $k_{4}=1$, corresponding to the quadruple $(3,1,5,1)$.

In both cases $(k_{1},k_{2},k_{3}) = (3,1,8)$ and $(k_{1},k_{2},k_{3}) = (4,1,5)$, we also have the single value $k_{4}=1$, which leads to the quadruples $(3,1,8,1)$ and $(4,1,5,1)$.

Therefore, the set $\mathcal{A}_{4}^{*}(3,1,0)$ consists of the quadruples
\begin{align}\label{lab_36}
& (k_{1},1,2,2),\quad k_{1}\geqslant 10;\notag\\
& (k_{1},1,2,4),\quad k_{1} = 6,7,8;\notag\\
& (2,1,k_{3},1),\quad k_{3}\equiv 2\pmod{3},\quad k_{3}\geqslant 8;\\
& (1,4,2,1);\notag\\
& (3,1,5,1),\;(3,1,8,1);\notag\\
& (4,1,5,1).\notag
\end{align}

In the \textsc{case} $r=5$, the set $\mathcal{A}_{5}^{\circ}(3,1,0)$ is contained of the tuples $(k_{1},\ldots,k_{5})$ such that $(k_{1},k_{2},k_{3},k_{4})$ belong to the list (\ref{lab_35}) and $k_{5}$ satisfies to the congruence
\[
\mathbb{K}_{4}(k_{2},k_{3},k_{4},k_{5}) = k_{5}\mathbb{K}_{3}(k_{2},k_{3},k_{4})-\mathbb{K}_{2}(k_{2},k_{3})\equiv 2\pmod{3}.
\]
However, in any case we have $\mathbb{K}_{2}(k_{2},k_{3})\equiv 1\pmod{3}$. Hence, the condition for $k_{5}$ takes the form
\begin{equation}\label{lab_37}
k_{5}\mathbb{K}_{3}(k_{2},k_{3},k_{4})\equiv 0\pmod{3}.
\end{equation}
In the case $(k_{1},k_{2},k_{3},k_{4}) = (k_{1},1,2,2)$, $k_{1}\geqslant 10$, one has $\mathbb{K}_{3}(k_{2},k_{3},k_{4})\equiv 1 \pmod{3}$. Thus, (\ref{lab_37}) implies that
$k_{5}\equiv 0\pmod{3}$. By lemma 2, the necessary condition of the non-emptiness of $\mathcal{T}(k_{1},\ldots,k_{5})$, that is, the non-emptiness of $\mathcal{T}(2,k_{5})$ implies that
$k_{5}\leqslant 5$. Therefore, we have the single possibility for $k_{5}$, namely, $k_{5} = 3$. Thus we get the family of tuples $(k_{1},1,2,2,3)$, where $k_{1}\geqslant 10$.

By lemma 4, both the regions $\mathcal{T}(k,1,2,2,3)$ and $\mathcal{T}(3,2,2,1,k)$ are simultaneously empty or non-empty.
By lemma 12, the region $\mathcal{T}(3,2,2,1,k)$ is non-empty only for $10\leqslant k\leqslant 16$. Hence, the set $\mathcal{A}_{5}^{\circ}(3,1,0)$ contains the tuples $(k_{1},1,2,2,3)$, $10\leqslant k_{1}\leqslant 16$.

All the remaining quadruples from the list (\ref{lab_36}) satisfy to the congruence $\mathbb{K}_{3}(k_{2},k_{3},k_{4})\equiv 0\pmod{3}$. Hence, the condition (\ref{lab_37}) is satisfied
by any $k_{5}\geqslant 1$.

If $(k_{1},k_{2},k_{3},k_{4}) = (k_{1},1,2,4)$, where $6\leqslant k_{1}\leqslant 8$, then the necessary condition, that is, the non-emptiness of the region $\mathcal{T}(2,4,k_{5})$, leads by lemma 3 (5b) to the single possible value $k_{5}=1$. The direct calculation shows that each of the regions $\mathcal{T}(k_{1},1,2,4)$, $6\leqslant k_{1}\leqslant 8$ is non-empty.

If $(k_{1},k_{2},k_{3},k_{4}) = (2,1,k_{3},1)$, where $k_{3}\equiv 2\pmod{3}$, $k_{3}\geqslant 8$, then, in the case $k_{3}\geqslant 4\cdot 5+2 =22$, by lemma 6, there is the single value  $k_{5} = 2$ corresponding to the non-empty region $\mathcal{T}(2,1,k_{3},1,k_{5})$. The direct calculation shows that the same holds for the values $k_{3} = 11, 14, 17, 20$. However, in the case $k_{3} = 8$ there exist two values of $k_{5}$ , namely, $k_{5} = 2$ and $k_{5} = 3$ corresponding to non-empty regions. Thus, the following tuples belong to the set $\mathcal{A}_{5}^{\circ}(3,1,0)$:
\begin{align*}
& (2,1,k_{3},1,2),\quad k_{3}\equiv 2\pmod{3},\quad k_{3}\geqslant 8;\\
& (2,1,8,1,3).
\end{align*}

If $(k_{1},k_{2},k_{3},k_{4}) = (1,4,2,1)$, then the non-emptiness of $\mathcal{T}(4,2,1,k_{5})$ is the necessary condition for the non-emptiness of the region $\mathcal{T}(k_{1},\ldots,k_{5}) = \mathcal{T}(1,4,2,1,k_{5})$.

By lemmas 4 and 8, the region $\mathcal{T}(4,2,1,k_{5})$ is non-empty only for $6\leqslant k_{5}\leqslant 8$. The non-emptiness of each region $\mathcal{T}(1,4,2,1,k)$, $6\leqslant k\leqslant 8$, is checked directly. Therefore, $\mathcal{A}_{5}^{\circ}(3,1,0)$ contains the tuples
\[
(1,4,2,1,6),\;(1,4,2,1,7),\;(1,4,2,1,8).
\]

If $(k_{1},k_{2},k_{3},k_{4}) = (3,1,5,1)$, then the necessary condition, that is, the non-emptiness of $\mathcal{T}(5,1,k_{5})$, leads by lemma 3 (2d) to the values $k_{5} = 3$ and $k_{5} = 4$. The calculation shows that the corresponding regions $\mathcal{T}(3,1,5,1,k_{5})$ are both non-empty. Hence, the set $\mathcal{A}_{5}^{\circ}(3,1,0)$ contains the tuples
\[
(3,1,5,1,3),\; (3,1,5,1,4).
\]
If $(k_{1},k_{2},k_{3},k_{4}) = (3,1,8,1)$, then the necessary condition, that is, the non-emptiness of $\mathcal{T}(8,1,k_{5})$, leads by lemma 3 (1),(2e) to the values $k_{5} = 2$ and $k_{5} = 3$. The direct calculation shows that the region $\mathcal{T}(3,1,8,1,k_{5})$ is non-empty only for $k_{5} = 2$. Hence, $\mathcal{A}_{5}^{\circ}(3,1,0)$ contains the tuple
\[
(3,1,8,1,2).
\]
Finally, if $(k_{1},k_{2},k_{3},k_{4}) = (4,1,5,1)$, then the necessary condition, that is, the non-emptiness of $\mathcal{T}(5,1,k_{5})$ leads by lemma 3 (2d) to the values $k_{5} = 3$, $k_{5}=4$. The direct calculation shows that the region $\mathcal{T}(4,1,5,1,k_{5})$ is non-empty only for $k_{5} = 3$. Therefore, the set $\mathcal{A}_{5}^{\circ}(3,1,0)$ contains the tuple
\[
(4,1,5,1,3).
\]
Thus, the set $\mathcal{A}_{5}^{\circ}(3,1,0)$ consists of the tuples
\begin{align}\label{lab_38}
& (k_{1},1,2,2,3),\quad 10\leqslant k_{1}\leqslant 16;\notag\\
& (k_{1},1,2,4,1),\quad k_{1} = 6,7,8;\notag\\
& (2,1,k_{3},1,2),\quad k_{3}\equiv 2\pmod{3},\quad k_{3}\geqslant 8;\notag\\
& (2,1,8,1,3);\\
& (1,4,2,1,k_{5}),\quad k_{5} = 6,7,8;\notag\\
& (3,1,5,1,3),\;(3,1,5,1,4),\; (3,1,8,1,2);\notag\\
& (4,1,5,1,3).\notag
\end{align}

At subsequent steps, in the structure of sets, <<stabilization>> occurs: for $\mathcal{A}_{r}^{*}(3,1,0)$ - starting with $r\geqslant 5$, for $\mathcal{A}_{r}^{\circ}(3,1,0)$ - starting with $r\geqslant 6$. In these cases, the description of the above sets looks quite simple. Namely, it appears that

(a) for $r\geqslant 5$, the set $\mathcal{A}_{r}^{*}(3,1,0)$ is formed by the tuples $(k_{1},1,2^{r-2})$, where $k_{1}\geqslant 4r-6$;

(b) $r\geqslant 6$, the set $\mathcal{A}_{r}^{\circ}(3,1,0)$ is formed by the tuples $(k_{1},1,2^{r-3},3)$, where $4r-10\leqslant k_{1}\leqslant 4r-4$.

We prove it by induction. Let $r=5$. In this case, the condition (\ref{lab_37}) is replaced by the following:
\begin{equation}\label{lab_39}
k_{5}\mathbb{K}_{3}(k_{2},k_{3},k_{4})\not\equiv 0\pmod{3}.
\end{equation}
As it was mentioned above, it fails for any $k_{5}$ for each quadruple from the list (\ref{lab_36}), excepting the tuples $(k_{1},k_{2},k_{3},k_{4}) = (k_{1},1,2,2)$, $k_{1}\geqslant 10$.
In the last case, (\ref{lab_39}) is reduced to $k_{5}\not\equiv 0\pmod{3}$.

The non-emptiness of the region $\mathcal{T}(2,2,k_{5})$ is the necessary condition for the non-emptiness of $\mathcal{T}(k_{1},\ldots,k_{5}) = \mathcal{T}(k_{1},1,2,2,k_{5})$.
By lemma 3 (3b), it holds only for $1\leqslant k_{5}\leqslant 3$. Excepting the value $k_{5} = 3$, we get $1\leqslant k_{5}\leqslant 2$. However, in the case $k_{5} = 1$, another necessary condition fails: the region $\mathcal{T}(1,2,2,k_{5}) = \mathcal{T}(1,2,2,1)$ is empty. Hence, $k_{5} = 2$ and the tuples under considering have the form $(k_{1},1,2,2,2)$, where $k_{1}\geqslant 10$. In view of lemmas 4 and 10, all the regions $\mathcal{T}(k_{1},1,2^{3})$, $k_{1}\geqslant 14$, are non-empty. This is the basis of induction.

Suppose now that $r\geqslant 6$. If the tuple $(k_{1},\ldots,k_{r})$ is contained in $\mathcal{A}_{r}^{*}(3,1,0)$, then, by induction, it has the form $(k_{1},1,2^{r-3},k_{r})$, where $k_{1}\geqslant 4r-10$. The condition (\ref{lab_25}) now reduces to the following:
\[
\mathbb{K}_{r-1}(1,2^{r-3},k_{r})\not\equiv 2\pmod{3}
\]
or, that is the same, to the condition
\begin{equation}\label{lab_40}
k_{r}\not\equiv 0\pmod{3}.
\end{equation}
The non-emptiness of the regions $\mathcal{T}(2,k_{r})$, $\mathcal{T}(2,2,k_{r})$  is the necessary condition for non-emptiness of $\mathcal{T}(k_{1},1,2^{r-3},k_{r})$. It narrows down the set of possible values of $k_{r}$ to the following: $1\leqslant k_{r}\leqslant 3$. Taking (\ref{lab_40}) into account, we conclude that $k_{r}=1$ or $k_{r} = 2$. However, in the case $k_{r}=1$ the necessary condition of non-emptiness of the region $\mathcal{T}(k_{1},1,2^{r-3},k_{r})$ fails because $\mathbb{K}_{r}(1,2^{r-3},k_{r})=0$.

Thus, we get the unique possibility $k_{r} = 2$ corresponding to the region $\mathcal{T}(k_{1},1,2^{r-2})$. Its non-emptiness for any $k_{1}\geqslant 4r-6$ follows now from lemma 10.
Hence, the set $\mathcal{A}_{r}^{*}(3,1,0)$ indeed consists of all tuples $(k_{1},1,2^{r-2})$, $k_{1}\geqslant 4r-6$. This finishes the description of the sets $\mathcal{A}_{r}^{*}(3,1,0)$ for $r\geqslant 5$.

Suppose now that the tuple $(k_{1},\ldots,k_{r})$ belongs to the set $\mathcal{A}_{r}^{\circ}(3,1,0)$, where $r\geqslant 6$. Then $(k_{1},\ldots,k_{r-1})\in \mathcal{A}_{r-1}^{*}(3,1,0)$ and,
in view of the above arguments,
\[
(k_{1},\ldots,k_{r-1}) = (k_{1},1,2^{r-3}),\quad\text{where}\quad k_{1}\geqslant 4(r-1)-6 = 4r-10.
\]
In this case, the condition (\ref{lab_10}) reduces to the following:
\[
\mathbb{K}_{r-1}(1,2^{r-3},k_{r})\equiv 2\pmod{3}\quad\text{or, that is the same,}\quad k_{r}\equiv 0\pmod{3}.
\]
The necessary condition of the region $\mathcal{T}(k_{1},1,2^{r-3},k_{r})$, that is, the non-emptiness of $\mathcal{T}(2,k_{r})$, yields $k_{r}\leqslant 5$, so we have $k_{r} = 3$.

As it was stated earlier, the regions $\mathcal{T}(3,2^{r-3},1,k_{r})$ are non-empty for $4r-10\leqslant k_{r}\leqslant 4r-4$. Hence, for $r\geqslant 6$, the set $\mathcal{A}_{r}^{\circ}(3,1,0)$ consists of the tuples $(k_{1},1,2^{r-3},3)$, $4r-10\leqslant k_{1}\leqslant 4r-4$. This finishes the description of the sets $\mathcal{A}_{r}^{\circ}(3,1,0)$.

We organize the above results as follows.
\vspace{0.3cm}

\textsc{Lemma 14.} \emph{The set} $\mathcal{A}_{r}^{\circ}(3,1,0)$ (\emph{and hence the set} $\mathcal{A}_{r}^{\circ}(3,2,0)$) \emph{consists of the following tuples}:

\begin{longtable}{|>{\fontsize{10}{9pt}\selectfont}c|>{\fontsize{10}{9pt}\selectfont}c|>{\fontsize{10}{9pt}\selectfont}c|}
\hline $r$ & $\mathbf{k} = (k_{1},\ldots,k_{r})$       & $\mathbb{K}_{r}(\mathbf{k})$  \\
\hline 2   & $(1,k_{2}), k_{2}\equiv 2\pmod{3}$, $k_{2}\geqslant 2$ & $k_{2}-1$\\
\hline 2   & $(2,2)$ & $3$\\
\hline 2   & $(3,2)$ & $5$\\
\hline 2   & $(4,2)$ & $7$\\
\hline 3   & $(1,k_{2},1), k_{2}\equiv 0\pmod{3}$, $k_{2}\geqslant 3$ & $k_{2}-2$\\
\hline 3   & $(2,1,k_{3})$, $k_{3}\equiv 0\pmod{3}$, $k_{3}\geqslant 6$ & $k_{3}-3$\\
\hline 3   & $(1,3,2),\;(2,3,1)$ & $3$\\
\hline 3   & $(2,3,2)$ & $8$\\
\hline 3   & $(3,1,6)$ & $9$\\
\hline 3   & $(k_{1},1,3)$, $4\leqslant k_{1}\leqslant 8$ & $2k_{1}-3$\\
\hline 4   & $(1,k_{2},1,2)$, $k_{2}\equiv 1\pmod{3}$, $k_{2}\geqslant 7$ & $k_{2}-3$\\
\hline 4   & $(2,1,k_{3},1)$, $k_{3}\equiv 1\pmod{3}$, $k_{3}\geqslant 7$ & $k_{3}-3$\\
\hline 4   & $(1,4,1,3),\;(3,1,4,1)$ & $3$\\
\hline 4   & $(1,4,1,4),\;(4,1,4,1)$ & $5$\\
\hline 4   & $(1,4,1,5),\;(5,1,4,1)$ & $7$\\
\hline 4   & $(2,4,1,3),\;(3,1,4,2)$ & $8$\\
\hline 4   & $(1,7,1,3),\;(3,1,7,1)$ & $9$\\
\hline 4   & $(2,4,1,4),\;(4,1,4,2)$ & $13$\\
\hline 4   & $(k_{1},1,2,3)$, $7\leqslant k_{1}\leqslant 12$ & $2k_{1}-5$\\
\hline 5   & $(2,1,k_{3},1,2)$, $k_{3}\equiv 2\pmod{3}$, $k_{2}\geqslant 8$ & $k_{3}-4$\\
\hline 5   & $(k_{1},1,2,2,3)$, $10\leqslant k_{1}\leqslant 16$ & $2k_{1}-7$\\
\hline 5   & $(k_{1},1,2,4,1)$, $6\leqslant k_{1}\leqslant 8$ & $2k_{1}-5$\\
\hline 5   & $(1,4,2,1,k_{5})$, $6\leqslant k_{5}\leqslant 8$ & $2k_{5}-5$\\
\hline 5   & $(3,1,5,1,3)$ & $8$\\
\hline 5   & $(2,1,8,1,3),\; (3,1,8,1,2)$ & $9$\\
\hline 5   & $(3,1,5,1,4),\;(4,1,5,1,3)$ & $13$\\
\hline $\geqslant 6$   & $(k_{1},1,2^{r-3},3)$, $4r-10\leqslant k_{1}\leqslant 4r-4$  & $2k_{1}-2r+3$\\
\hline
\end{longtable}
\vspace{0.5cm}

\section{Proof of the main assertion}
\vspace{0.5cm}

By lemma 1, in each case $c_{0} = 1$, $c_{0} = 2$, for any $r\geqslant 1$ we have:
\begin{multline*}
\nu(Q;r,3,c_{0}) = \nu(r,3,c_{0}) + O\biggl(\frac{\ln{Q}}{Q}\biggr),\quad \nu(r,3,c_{0}) = \frac{2}{3}\,\mathfrak{c}_{r},\\
\mathfrak{c}_{r} = \biggl(\,\sum\limits_{\mathbf{k}\in\,\mathcal{A}_{r}^{\circ}(3,1,0)} + \sum\limits_{\mathbf{k}\in\,\mathcal{A}_{r}^{\circ}(3,1,2)}\;\biggr)|\mathcal{T}(\mathbf{k})|
= \biggl(\,\sum\limits_{\mathbf{k}\in\,\mathcal{A}_{r}^{\circ}(3,2,0)} + \sum\limits_{\mathbf{k}\in\,\mathcal{A}_{r}^{\circ}(3,2,1)}\;\biggr)|\mathcal{T}(\mathbf{k})|.
\end{multline*}
\textsc{Case} $r=1$. By lemma 13,
\begin{multline*}
\nu(1,3,1) = \frac{2}{3}\,\sum\limits_{\mathbf{k}\in\,\mathcal{A}_{1}^{\circ}(3,1,2)}|\mathcal{T}(\mathbf{k})| = \frac{2}{3}\,\sum\limits_{\substack{k\equiv 1\;(\mmod{3}) \\ k\geqslant 1}}|\mathcal{T}(k)| = \frac{2}{3}\biggl(|\mathcal{T}(1)| + \sum\limits_{k=1}^{+\infty}|\mathcal{T}(3k+1)|\biggr) =\\
= \frac{2}{3}\biggl(\frac{1}{6} + \sum\limits_{k=1}^{+\infty}\frac{4}{(3k+1)(3k+2)(3k+3)}\biggr) = \frac{2\pi}{3\sqrt{3}} - \frac{2}{3}\ln{3} - \frac{1}{3} = 0.14345\,80504\ldots
\end{multline*}
\textsc{Case} $r=2$. By lemmas 13 and 14,
\begin{multline*}
\nu(2,3,1) = \frac{2}{3}\biggl\{\,\sum\limits_{\mathbf{k}\in\,\mathcal{A}_{2}^{\circ}(3,1,2)}|\mathcal{T}(\mathbf{k})| +  \,\sum\limits_{\mathbf{k}\in\,\mathcal{A}_{2}^{\circ}(3,1,0)}|\mathcal{T}(\mathbf{k})|\biggr\}= \\
= \frac{2}{3}\biggl\{\,\sum\limits_{\substack{k_{1}\equiv 2\;(\mmod{3}) \\ k_{1}\geqslant 2}}|\mathcal{T}(k_{1},1)| + |\mathcal{T}(2,2)|  + |\mathcal{T}(2,3)|  + |\mathcal{T}(2,4)| + \\
+ \,\sum\limits_{\substack{k_{2}\equiv 2\;(\mmod{3}) \\ k_{2}\geqslant 2}}|\mathcal{T}(1,k_{2})|  + |\mathcal{T}(2,2)| + |\mathcal{T}(3,2)|  + |\mathcal{T}(4,2)|\biggr\} =\\
= \frac{4}{3}\biggl\{\,\sum\limits_{\substack{k_{1}\equiv 2\;(\mmod{3}) \\ k_{1}\geqslant 2}}|\mathcal{T}(k_{1},1)| + \frac{1}{10}  + \frac{1}{35} + \frac{1}{210}\biggr\}.
\end{multline*}
By lemma 7 (a), $|\mathcal{T}(k_{1},1)| = |\mathcal{T}(k_{1})|$ for $k_{1}\geqslant 5$. Therefore,
\begin{multline*}
\nu(2,3,1) = \frac{4}{3}\biggl\{\,\sum\limits_{\substack{k_{1}\equiv 2\;(\mmod{3}) \\ k_{1}\geqslant 5}}|\mathcal{T}(k_{1})| + |\mathcal{T}(2,1)| + \frac{2}{15}\biggr\} = \\
= \frac{4}{3}\biggl\{\;\sum\limits_{k=1}^{+\infty}\frac{4}{(3k+2)(3k+3)(3k+4)}+\frac{1}{6}\biggr\} = \frac{8}{3}(\ln{3}-1) = 0.26296\,62031\ldots
\end{multline*}
\textsc{Case} $r=3$. We have
\[
\nu(3,3,1) = \frac{2}{3}\biggl\{\,\sum\limits_{\mathbf{k}\in\,\mathcal{A}_{3}^{\circ}(3,1,2)}|\mathcal{T}(\mathbf{k})| +  \,\sum\limits_{\mathbf{k}\in\,\mathcal{A}_{3}^{\circ}(3,1,0)}|\mathcal{T}(\mathbf{k})|\biggr\} = \frac{2}{3}(\sigma_{1}+\sigma_{2}),
\]
where the notations are obvious. By lemma 13,
\[
\sigma_{2} = \sum\limits_{\substack{k_{1}\equiv 0\;(\mmod{3}) \\ k_{1}\geqslant 3}}|\mathcal{T}(k_{1},1,2)| + \sum\limits_{k_{3}=4}^{8}|\mathcal{T}(3,1,k_{3})| + |\mathcal{T}(6,1,3)|.
\]
By lemma 7 (b), $|\mathcal{T}(k_{1},1,2)| = |\mathcal{T}(k_{1})|$ for $k_{1}\geqslant 9$. Therefore,
\begin{multline*}
\sigma_{1} = \sum\limits_{\substack{k_{1}\equiv 0\;(\mmod{3}) \\ k_{1}\geqslant 9}}|\mathcal{T}(k_{1})| + |\mathcal{T}(3,1,2)| + |\mathcal{T}(6,1,2)| + \frac{221}{4620} = \\
=\sum\limits_{k=3}^{+\infty}\frac{4}{3k(3k+1)(3k+2)} + \frac{1}{20} = \frac{104}{35}-\frac{\pi}{\sqrt{3}}-\ln{3}.
\end{multline*}
Similarly, by lemma 14 we get
\begin{multline*}
\sigma_{2} = \sum\limits_{\substack{k_{2}\equiv 0\;(\mmod{3}) \\ k_{2}\geqslant 3}}|\mathcal{T}(1,k_{2},1)| + \sum\limits_{\substack{k_{3}\equiv 0\;(\mmod{3}) \\ k_{3}\geqslant 6}}|\mathcal{T}(2,1,k_{3})| + |\mathcal{T}(1,3,2)| + |\mathcal{T}(2,3,1)| + \\
+ |\mathcal{T}(2,3,2)| + |\mathcal{T}(3,1,6)| + \sum\limits_{k_{1}=4}^{8}|\mathcal{T}(k_{1},1,3)|.
\end{multline*}
In view of lemma 7 (f),(c), $|\mathcal{T}(1,k,1)| = |\mathcal{T}(k)|$ for $k\geqslant 6$,
$|\mathcal{T}(2,1,k)| = |\mathcal{T}(k,1,2)| = |\mathcal{T}(k)|$ for $k\geqslant 9$, and hence
\begin{multline*}
\sigma_{2} = \sum\limits_{\substack{k_{2}\equiv 0\;(\mmod{3}) \\ k_{2}\geqslant 6}}|\mathcal{T}(k_{2})| + \sum\limits_{\substack{k_{3}\equiv 0\;(\mmod{3}) \\ k_{3}\geqslant 9}}|\mathcal{T}(k_{3})| + |\mathcal{T}(1,3,1)| + |\mathcal{T}(2,1,3)|  + |\mathcal{T}(2,1,6)| + \frac{113}{1155} = \\
= \sum\limits_{k=2}^{+\infty}\frac{4}{3k(3k+1)(3k+2)} + \sum\limits_{k=3}^{+\infty}\frac{4}{3k(3k+1)(3k+2)} + \frac{7}{60} = \frac{209}{35} - \frac{2\pi}{\sqrt{3}} - 2\ln{3}.
\end{multline*}
Thus we finally get
\[
\nu(3,3,1) = \frac{626}{105} - \frac{2\pi}{\sqrt{3}} - 2\ln{3} = 0.13708\,14561\ldots.
\]
\textsc{Case} $r=4$. We have
\[
\nu(4,3,1) = \frac{2}{3}(\sigma_{1}+\sigma_{2}),\quad \sigma_{1} = \sum\limits_{\mathbf{k}\in\,\mathcal{A}_{4}^{\circ}(3,1,2)}|\mathcal{T}(\mathbf{k})|,
\quad \sigma_{2} = \sum\limits_{\mathbf{k}\in\,\mathcal{A}_{4}^{\circ}(3,1,0)}|\mathcal{T}(\mathbf{k})|.
\]
By lemma 13,
\[
\sigma_{1} = \sum\limits_{k_{4}=7}^{12}|\mathcal{T}(3,2,1,k_{4})| = \frac{1}{70}.
\]
Further, by lemma 14,
\begin{multline*}
\sigma_{2} = \sum\limits_{\substack{k_{2}\equiv 1\;(\mmod{3}) \\ k_{2}\geqslant 7}}|\mathcal{T}(1,k_{2},1,2)| + \sum\limits_{\substack{k_{3}\equiv 1\;(\mmod{3}) \\ k_{3}\geqslant 7}}|\mathcal{T}(2,1,k_{3},1)| + \\
+ |\mathcal{T}(1,4,1,3)| + |\mathcal{T}(1,4,1,4)|  + |\mathcal{T}(1,4,1,5)| + |\mathcal{T}(1,7,1,3)| + \\
+ |\mathcal{T}(2,4,1,3)| + |\mathcal{T}(2,4,1,4)|  + |\mathcal{T}(3,1,4,1)| + |\mathcal{T}(3,1,4,2)| + \\
+ |\mathcal{T}(3,1,7,1)| + |\mathcal{T}(4,1,4,1)|  + |\mathcal{T}(4,1,4,2)| + |\mathcal{T}(5,1,4,1)| + \\
+ \sum\limits_{k_{1}=7}^{12}|\mathcal{T}(k_{1},1,2,3)| = \sum\limits_{\substack{k_{2}\equiv 1\;(\mmod{3}) \\ k_{1}\geqslant 7}}|\mathcal{T}(1,k_{2},1,2)| + \frac{710}{9009}.
\end{multline*}
Using lemma 4 and lemma 7 (d), we conclude that $|\mathcal{T}(1,k,1,2)| = |\mathcal{T}(2,1,k,1)| = |\mathcal{T}(k)|$ for $k\geqslant 9$. Hence,
\[
\sigma_{2} = \sum\limits_{k=3}^{+\infty}\frac{4}{(3k+1)(3k+2)(3k+3)} + \frac{17}{2002} + \frac{710}{9009} = \frac{2\pi}{\sqrt{3}} - \frac{93}{70} - 2\ln{3}.
\]
Finally we obtain
\[
\nu(4,3,1) = \frac{4\pi}{3\sqrt{3}} - \frac{92}{105} - \frac{4}{3}\ln{3} = 0.07739\,22912\ldots.
\]
\textsc{Case} $r=5$. Here we have
\[
\nu(5,3,1) = \frac{2}{3}(\sigma_{1}+\sigma_{2}),\quad \sigma_{1} = \sum\limits_{\mathbf{k}\in\,\mathcal{A}_{5}^{\circ}(3,1,2)}|\mathcal{T}(\mathbf{k})|,
\quad \sigma_{2} = \sum\limits_{\mathbf{k}\in\,\mathcal{A}_{5}^{\circ}(3,1,0)}|\mathcal{T}(\mathbf{k})|.
\]
By lemma 13,
\[
\sigma_{1} = \sum\limits_{k_{5}=10}^{16}|\mathcal{T}(3,2,2,1,k_{5})| = \frac{2}{315}.
\]
Further, lemma 14 implies that
\begin{multline*}
\sigma_{2} = \sum\limits_{\substack{k_{3}\equiv 2\;(\mmod{3}) \\ k_{3}\geqslant 7}}|\mathcal{T}(2,1,k_{3},1,2)| + \sum\limits_{k_{1}=10}^{16}|\mathcal{T}(k_{1},1,2,2,3)| +
\sum\limits_{k_{1}=6}^{8}|\mathcal{T}(k_{1},1,2,4,1)| + \\
+ |\mathcal{T}(2,1,8,1,3)| + \sum\limits_{k_{5}=6}^{8}|\mathcal{T}(1,4,2,1,k_{5})| + \\
+ |\mathcal{T}(3,1,5,1,3)|  + |\mathcal{T}(3,1,5,1,4)| + |\mathcal{T}(3,1,8,1,2)| + |\mathcal{T}(4,1,5,1,3)| = \\
 = \sum\limits_{\substack{k_{3}\equiv 2\;(\mmod{3}) \\ k_{3}\geqslant 8}}|\mathcal{T}(2,1,k_{3},1,2)| + \frac{10}{273}.
\end{multline*}
In view of lemma 7 (e), one has $|\mathcal{T}(2,1,k_{3},1,2)| = |\mathcal{T}(k_{3})|$ for $k_{3}\geqslant 9$. Therefore,
\begin{multline*}
\sigma_{2} = \sum\limits_{\substack{k_{3}\equiv 2\;(\mmod{3}) \\ k_{3}\geqslant 11}}|\mathcal{T}(k_{3})|  + \frac{1}{260}+ \frac{10}{273} =\\
= \sum\limits_{k=3}^{+\infty}\frac{4}{(3k+2)(3k+3)(3k+4)} + \frac{17}{420} = 2\ln{3} - \frac{271}{126}.
\end{multline*}
Hence,
\[
\nu(5,3,1) = \frac{4}{3}\ln{3} - \frac{193}{135} = 0.03518\,67552\ldots
\]
\textsc{Case} $r\geqslant 6$. Lemmas 13, 14 together with lemma 4 imply that
\[
\nu(r,3,1) = \frac{4}{3}\,\sigma,\quad \sigma = \sum\limits_{k_{1}=4r-10}^{4r-4}|\mathcal{T}(k_{1},1,2^{r-3},3)| = \sum\limits_{k_{r}=4r-10}^{4r-4}|\mathcal{T}(3,2^{r-3},1,k_{r})|.
\]
By lemma 12, the area of the union of the regions $\mathcal{T}(3,2^{r-3},1,k_{r})$, $4r-10\leqslant k_{r}\leqslant 4r-4$ equals to
\[
\frac{2}{(2r-5)(2r-3)(2r-1)}.
\]
Therefore,
\[
\nu(r,3,1) = \frac{8}{3(2r-5)(2r-3)(2r-1)}.
\]
Finally we note that
\[
\sum\limits_{r=6}^{+\infty}\frac{8}{3(2r-5)(2r-3)(2r-1)} = \frac{2}{189}.
\]
Hence, summing all the above proportions, we get
\begin{multline*}
\sum\limits_{r=1}^{+\infty}\nu(r,3,1) = \biggl(\frac{2\pi}{3\sqrt{3}} - \frac{2}{3}\ln{3} - \frac{1}{3}\biggr) +
\frac{8}{3}(\ln{3}-1) + \biggl(\frac{626}{105} - \frac{2\pi}{\sqrt{3}} - 2\ln{3}\biggr) + \\
+ \biggl(\frac{4\pi}{3\sqrt{3}} - \frac{92}{105} - \frac{4}{3}\ln{3}\biggr) + \biggl(\frac{4}{3}\ln{3} - \frac{193}{135}\biggr) + \frac{2}{189} = \frac{124}{189} + \frac{2}{189} = \frac{2}{3}.
\end{multline*}
Thus,
\[
\nu(0,3,1) = 1-\frac{2}{3} = \frac{1}{3}.
\]
Theorem is proved.
\vspace{1cm}

\noindent
\large
\textbf{Appendix I. Precise description of the polygons $\boldsymbol{\mathcal{T}(3,2,1,k)}$, $\boldsymbol{7\leqslant k\leqslant 12}$}
\vspace{0.5cm}

\normalsize
\textsc{Case} $k=7$: the region $\mathcal{T}(3,2,1,7)$ is described by the system of the inequalities
\begin{equation*}
\begin{cases}
\displaystyle \frac{4}{7} < x\leqslant \frac{10}{17},\quad 1-x<y\leqslant \frac{1+2x}{5},\\[12pt]
\displaystyle \frac{10}{17} < x\leqslant \frac{3}{4},\quad \frac{1+6x}{11}<y\leqslant \frac{1+2x}{5}
\end{cases}
\end{equation*}
and is the triangle with the vertices
\[
A=\biggl(\frac{4}{7},\frac{3}{7}\biggr),\quad B=\biggl(\frac{3}{4},\frac{1}{2}\biggr),\quad C=\biggl(\frac{10}{17},\frac{7}{17}\biggr) \quad\text{and with area}\quad \frac{1}{476}.
\]
\vspace{0.3cm}

\textsc{Case} $k=8$: the region $\mathcal{T}(3,2,1,8)$ is described by the system of the inequalities
\begin{equation*}
\begin{cases}
\displaystyle \frac{10}{17} < x\leqslant \frac{3}{5},\quad 1-x<y\leqslant \frac{1+6x}{11},\\[12pt]
\displaystyle \frac{3}{5} < x\leqslant \frac{3}{4},\quad \frac{1+7x}{13}<y\leqslant \frac{1+6x}{11},\\[12pt]
\displaystyle \frac{3}{4} < x\leqslant \frac{8}{9},\quad \frac{1+7x}{13}<y\leqslant \frac{1+2x}{5},
\end{cases}
\end{equation*}
and is the quadrangle with the vertices
\[
A=\biggl(\frac{10}{17},\frac{7}{17}\biggr),\quad B=\biggl(\frac{3}{4},\frac{1}{2}\biggr),\quad C=\biggl(\frac{8}{9},\frac{5}{9}\biggr), \quad D=\biggl(\frac{3}{5},\frac{2}{5}\biggr)
\quad\text{and with area}\quad \frac{13}{3\,060}.
\]
\vspace{0.3cm}

\textsc{Case} $k=9$: the region $\mathcal{T}(3,2,1,9)$ is described by the system of the inequalities
\begin{equation*}
\begin{cases}
\displaystyle \frac{3}{5} < x\leqslant \frac{8}{11},\quad \frac{1+3x}{7}<y\leqslant \frac{1+7x}{13},\\[12pt]
\displaystyle \frac{8}{11} < x\leqslant \frac{8}{9},\quad \frac{1+8x}{15}<y\leqslant \frac{1+7x}{13},\\[12pt]
\displaystyle \frac{8}{9} < x\leqslant 1,\quad \frac{1+8x}{15}<y\leqslant \frac{1+2x}{5},
\end{cases}
\end{equation*}
and is the quadrangle with the vertices
\[
A=\biggl(\frac{3}{5},\frac{2}{5}\biggr),\quad B=\biggl(\frac{8}{9},\frac{5}{9}\biggr),\quad C=\biggl(1,\frac{3}{5}\biggr), \quad D=\biggl(\frac{8}{11},\frac{5}{11}\biggr)
\quad\text{and with area}\quad \frac{2}{495}.
\]
\vspace{0.3cm}

\textsc{Case} $k=10$: the region $\mathcal{T}(3,2,1,10)$ is described by the system of the inequalities
\begin{equation*}
\begin{cases}
\displaystyle \frac{8}{11} < x\leqslant \frac{5}{6},\quad \frac{1+3x}{7}<y\leqslant \frac{1+8x}{15},\\[12pt]
\displaystyle \frac{5}{6} < x\leqslant 1,\quad \frac{1+9x}{17}<y\leqslant \frac{1+8x}{15},
\end{cases}
\end{equation*}
and is the quadrangle with the vertices
\[
A=\biggl(\frac{8}{11},\frac{5}{11}\biggr),\quad B=\biggl(1,\frac{3}{5}\biggr),\quad C=\biggl(1,\frac{10}{17}\biggr), \quad D=\biggl(\frac{5}{6},\frac{1}{2}\biggr)
\quad\text{and with area}\quad \frac{7}{2\,805}.
\]
\vspace{0.3cm}

\textsc{Case} $k=11$: the region $\mathcal{T}(3,2,1,11)$ is described by the system of the inequalities
\begin{equation*}
\begin{cases}
\displaystyle \frac{5}{6} < x\leqslant \frac{12}{13},\quad \frac{1+3x}{7}<y\leqslant \frac{1+9x}{17},\\[12pt]
\displaystyle \frac{12}{13} < x\leqslant 1,\quad \frac{1+10x}{19}<y\leqslant \frac{1+9x}{17},
\end{cases}
\end{equation*}
and is the quadrangle with the vertices
\[
A=\biggl(\frac{5}{6},\frac{1}{2}\biggr),\quad B=\biggl(1,\frac{10}{17}\biggr),\quad C=\biggl(1,\frac{11}{19}\biggr), \quad D=\biggl(\frac{12}{13},\frac{7}{13}\biggr)
\quad\text{and with area}\quad \frac{14}{12\,597}.
\]
\vspace{0.3cm}

\textsc{Case} $k=12$: the region $\mathcal{T}(3,2,1,12)$ is described by the system of the inequalities
\[
\frac{12}{13} < x\leqslant 1,\quad \frac{1+3x}{7}<y\leqslant \frac{1+10x}{19}
\]
and is the triangle with the vertices
\[
A=\biggl(\frac{12}{13},\frac{7}{13}\biggr),\quad B=\biggl(1,\frac{11}{19}\biggr),\quad C=\biggl(1,\frac{4}{7}\biggr)
\quad\text{and with area}\quad \frac{1}{3\,458}.
\]
\vspace{0.5cm}

\noindent
\large \textbf{Appendix II. Precise description of the polygons} \\ $\boldsymbol{\mathcal{T}(3,2^{r-3},1,k)}$, $\boldsymbol{4r-10\leqslant k\leqslant 4r-4}$
\vspace{0.5cm}

\normalsize
\textsc{Case} $k=4r-10$: the region  $\mathcal{T}(3,2^{r-3},1,4r-10)$ is described by the system of the inequalities
\begin{equation*}
\begin{cases}
\displaystyle \frac{2r-7}{2r-5} < x\leqslant \frac{6r-19}{6r-13},\quad \frac{1+x}{4}<y\leqslant \frac{1+(r-2)x}{2r-3},\\[12pt]
\displaystyle \frac{6r-19}{6r-13} < x\leqslant \frac{4r-12}{4r-9},\quad \frac{1+(3r-7)x}{6r-15}<y\leqslant \frac{1+(r-2)x}{2r-3}
\end{cases}
\end{equation*}
and is the triangle with the vertices
\begin{multline*}
A = \biggl(\frac{2r-7}{2r-5},\frac{r-3}{2r-5}\biggr),\quad B = \biggl(\frac{4r-12}{4r-9},\frac{2r-5}{4r-9}\biggr),\quad C = \biggl(\frac{6r-19}{6r-13},\frac{3r-8}{6r-13}\biggr)\\
\text{and with area}\quad \frac{1}{2(2r-5)(4r-9)(6r-13)}.
\end{multline*}

\vspace{0.3cm}

\textsc{Case} $k=4r-9$: the region $\mathcal{T}(3,2^{r-3},1,4r-9)$ is described by the system of the inequalities
\begin{equation*}
\begin{cases}
\displaystyle \frac{6r-19}{6r-13} < x\leqslant \frac{6r-17}{6r-11},\quad \frac{1+x}{4}<y\leqslant \frac{1+(3r-7)x}{6r-15} \\[12pt]
\displaystyle \frac{6r-17}{6r-11} < x\leqslant \frac{4r-12}{4r-9},\quad \frac{1+(3r-6)x}{6r-13}<y\leqslant \frac{1+(3r-7)x}{6r-15} \\[12pt]
\displaystyle \frac{4r-12}{4r-9} < x\leqslant \frac{2r-5}{2r-4},\quad \frac{1+(3r-6)x}{6r-13}<y\leqslant \frac{1+(r-2)x}{2r-3}
\end{cases}
\end{equation*}
and is the quadrangle with the vertices
\begin{multline*}
A = \biggl(\frac{6r-19}{6r-13},\frac{3r-8}{6r-13}\biggr),\quad B = \biggl(\frac{4r-12}{4r-9},\frac{2r-5}{4r-9}\biggr),\quad C = \biggl(\frac{2r-5}{2r-4},\frac{1}{2}\biggr),\\
D = \biggl(\frac{6r-17}{6r-11},\frac{3r-7}{6r-11}\biggr)\quad \text{and with area}\quad \frac{9r-19}{2(r-2)(4r-9)(6r-13)(6r-11)}.
\end{multline*}
\vspace{0.3cm}

\textsc{Case} $k=4r-8$: the region $\mathcal{T}(3,2^{r-3},1,4r-8)$ is described by the system of the inequalities
\begin{equation*}
\begin{cases}
\displaystyle \frac{6r-17}{6r-11} < x\leqslant \frac{2r-5}{2r-3},\quad \frac{1+x}{4}<y\leqslant \frac{1+(3r-6)x}{6r-13} \\[12pt]
\displaystyle \frac{2r-5}{2r-3} < x\leqslant \frac{2r-5}{2r-4},\quad \frac{1+(3r-5)x}{6r-11}<y\leqslant \frac{1+(3r-6)x}{6r-13} \\[12pt]
\displaystyle \frac{2r-5}{2r-4} < x\leqslant \frac{4r-8}{4r-7},\quad \frac{1+(3r-5)x}{6r-11}<y\leqslant \frac{1+(r-2)x}{2r-3}
\end{cases}
\end{equation*}
and is the quadrangle with the vertices
\begin{multline*}
A = \biggl(\frac{6r-17}{6r-11},\frac{3r-7}{6r-11}\biggr),\quad B = \biggl(\frac{2r-5}{2r-4},\frac{1}{2}\biggr),\quad C = \biggl(\frac{4r-8}{4r-7},\frac{2r-3}{4r-7}\biggr),\\
D = \biggl(\frac{2r-5}{2r-3},\frac{r-2}{2r-3}\biggr)\quad \text{and with area}\quad \frac{5r-9}{2(r-2)(2r-3)(4r-7)(6r-11)}.
\end{multline*}
\vspace{0.3cm}

\textsc{Case} $k=4r-7$: the region $\mathcal{T}(3,2^{r-3},1,4r-7)$ is described by the system of the inequalities
\begin{equation*}
\begin{cases}
\displaystyle \frac{2r-5}{2r-3} < x\leqslant \frac{4r-8}{4r-7},\quad \frac{1+(r-1)x}{2r-1}<y\leqslant \frac{1+(3r-5)x}{6r-11} \\[12pt]
\displaystyle \frac{4r-8}{4r-5} < x\leqslant \frac{4r-8}{4r-7},\quad \frac{1+(3r-4)x}{6r-9}<y\leqslant \frac{1+(3r-5)x}{6r-11} \\[12pt]
\displaystyle \frac{4r-8}{4r-7} < x\leqslant 1,\quad \frac{1+(3r-4)x}{6r-9}<y\leqslant \frac{1+(r-2)x}{2r-3}
\end{cases}
\end{equation*}
and is the quadrangle with the vertices
\begin{multline*}
A = \biggl(\frac{2r-5}{2r-3},\frac{r-2}{2r-3}\biggr),\quad B = \biggl(\frac{4r-8}{4r-7},\frac{2r-3}{4r-7}\biggr),\quad C = \biggl(1,\frac{r-1}{2r-3}\biggr),\\
D = \biggl(\frac{4r-8}{4r-5},\frac{2r-3}{4r-5}\biggr)\quad \text{and with area}\quad \frac{2}{(2r-3)(4r-7)(4r-5)}.
\end{multline*}
\vspace{0.3cm}

\textsc{Case} $k=4r-6$: the region $\mathcal{T}(3,2^{r-3},1,4r-6)$ is described by the system of the inequalities
\begin{equation*}
\begin{cases}
\displaystyle \frac{4r-8}{4r-5} < x\leqslant \frac{2r-3}{2r-2},\quad \frac{1+(r-1)x}{2r-1}<y\leqslant \frac{1+(3r-4)x}{6r-9} \\[12pt]
\displaystyle \frac{2r-3}{2r-2} < x\leqslant 1,\quad \frac{1+(3r-3)x}{6r-7}<y\leqslant \frac{1+(3r-4)x}{6r-9}
\end{cases}
\end{equation*}
and is the quadrangle with the vertices
\begin{multline*}
A = \biggl(\frac{4r-8}{4r-5},\frac{2r-3}{4r-5}\biggr),\quad B = \biggl(1,\frac{r-1}{2r-3}\biggr),\quad C = \biggl(1,\frac{3r-2}{6r-7}\biggr),\\
D = \biggl(\frac{2r-3}{2r-2},\frac{1}{2}\biggr)\quad \text{and with area}\quad \frac{5r-6}{2(r-1)(2r-3)(4r-5)(6r-7)}.
\end{multline*}
\vspace{0.3cm}

\textsc{Case} $k=4r-5$: the region $\mathcal{T}(3,2^{r-3},1,4r-5)$ is described by the system of the inequalities
\begin{equation*}
\begin{cases}
\displaystyle \frac{2r-3}{2r-2} < x\leqslant \frac{4r-4}{4r-3},\quad \frac{1+(r-1)x}{2r-1}<y\leqslant \frac{1+(3r-3)x}{6r-7} \\[12pt]
\displaystyle \frac{4r-4}{4r-3} < x\leqslant 1,\quad \frac{1+(3r-2)x}{6r-5}<y\leqslant \frac{1+(3r-3)x}{6r-7}
\end{cases}
\end{equation*}
and is the quadrangle with the vertices
\begin{multline*}
A = \biggl(\frac{2r-3}{2r-2},\frac{1}{2}\biggr),\quad B = \biggl(1,\frac{3r-2}{6r-7}\biggr),\quad C = \biggl(1,\frac{3r-1}{6r-5}\biggr),\\
D = \biggl(\frac{4r-4}{4r-3},\frac{2r-1}{4r-3}\biggr)\quad \text{and with area}\quad \frac{9r-8}{2(r-1)(4r-3)(6r-7)(6r-5)}.
\end{multline*}
\vspace{0.3cm}

\textsc{Case} $k=4r-4$: the region $\mathcal{T}(3,2^{r-3},1,4r-4)$ is described by the system of the inequalities
\begin{equation*}
\displaystyle \frac{4r-4}{4r-3} < x\leqslant 1,\quad \frac{1+(r-1)x}{2r-1}<y\leqslant \frac{1+(3r-2)x}{6r-5}
\end{equation*}
and is the triangle with the vertices
\begin{multline*}
A = \biggl(\frac{4r-4}{4r-3},\frac{2r-1}{4r-3}\biggr),\quad B = \biggl(1,\frac{3r-1}{6r-5}\biggr),\quad C = \biggl(1,\frac{r}{2r-1}\biggr),\\
\text{and with area}\quad \frac{1}{2(2r-1)(4r-3)(6r-5)}.
\end{multline*}

\renewcommand{\refname}{\Large{References}}

\end{document}